\documentclass[11pt]{amsart}
\usepackage{amsmath}

\usepackage{amssymb}

\newtheorem{theorem}{Theorem}[section]

\newtheorem{lemma}[theorem]{Lemma}

\newtheorem{corollary}[theorem]{Corollary}

\theoremstyle{definition}
\newtheorem{definition}[theorem]{Definition}

\theoremstyle{remark}
\newtheorem{remark}[theorem]{Remark}
\newtheorem{conclusion}[theorem]{Conclusion}
\newtheorem{notation}[theorem]{Notation}

\newtheorem{conjecture}[theorem]{Conjecture}

\newcount\skewfactor
\def\mathunderaccent#1#2 {\let\theaccent#1\skewfactor#2
\mathpalette\putaccentunder}
\def\putaccentunder#1#2{\oalign{$#1#2$\crcr\hidewidth
\vbox to.2ex{\hbox{$#1\skew\skewfactor\theaccent{}$}\vss}\hidewidth}}

\def\smallbox#1{\leavevmode\thinspace\hbox{\vrule\vtop{\vbox
   {\hrule\kern1pt\hbox{\vphantom{\tt/}\thinspace{\tt#1}\thinspace}}
   \kern1pt\hrule}\vrule}\thinspace}

\newcommand\includedin\subseteq
\newcommand\intersect\cap
\newcommand\union\cup

\newcommand{\cf}{{\rm cf}}

\newcommand{\st}{{such that}}
\newcommand{\sq}{{sequence}}

\newcommand{\rest}{{\restriction}}

\newcommand{\real}{{R}}
\newcommand{\rationals}{{Q}}
\newcommand{\minus}{{--}}
\newcommand{\sminus}{{\setminus}}

\newcommand{\conc}{{}^\frown\!}

\usepackage[hidelinks]{hyperref}
\begin{document}
\makeatletter\def\shfiuwefootnote{\gdef\@thefnmark{}\@footnotetext}\makeatother\shfiuwefootnote{Version 2023-05-01. See \url{https://shelah.logic.at/papers/42/} for possible updates.}

\title{The Monadic Theory of Order}
% Sh:54

\author{Saharon Shelah}
\address{Institute of Mathematics
 The Hebrew University of Jerusalem
 Jerusalem 91904, Israel
 and  Department of Mathematics
 Rutgers University
 New Brunswick, NJ 08854, USA}
\email{shelah@math.huji.ac.il}
% \urladdr{http://www.math.rutgers.edu/\char`\~shelah} 
%%	Notimplemented locally
\thanks{First typed: September 1975. \\
Research supported by the United States-Israel Binational
Science Foundation. Publication sh:42}

%\subjclass{}
\keywords{monadic order, decidability}
\date{2021-07-08}
% 2021-07-08 09:12 minor araise date
% 2021-06-30 09:45 --% 2021-06-30 10:54 tiqunet hapnayot-% 2021-06-30 10:55 
% 2021-06-15 06:29 -2021-06-15 08:23 
% 2021-06-14 07:53 -- % 2021-06-14 18:52 - hapsaqot 5abalharbeh
% 2021-06-14 06:23 --2021-06-14 07:18 xon - naseq ++
% 2021-06-14 02:31 --% 2021-06-14 03:40 vad vamud 11

\begin{abstract}
%%% Put abstract here:
We deal with the monadic (second-order) theory of order. We prove all 
known results in a unified way, show a general way of reduction, prove more 
results and show the limitation on extending them. We prove (CH) that the monadic theory 
of the real order is undecidable. Our methods are model-theoretic, and we do not use 
automaton theory. 

This is a slightly corrected version of a very old work.
\end{abstract}
\maketitle

\section{Introduction}

% 2021-06-14 02:32 \smallbox{ siman h-diAz hmopiva val makaj hasipra 3 lo jamij. simanty b ?}

The monadic logic is first order logic when we add variables ranging over sets, and 
allow quantification over them. If pairing functions are available this is essentially 
second order logic. The monadic theory of a class $K$ of $L$-models is 
$\{\psi:\psi$ is a sentence in monadic logic, satisfied by any member of $K\}$.

Here we shall investigate cases where the members of $K$ are linear orders 
(with one-place predicates). 
    % 2021-06-14 02:33 \smallbox{kAn?}
\[{\rm Contents}\]
Section
\begin{enumerate}
\item[(0)] Introduction ......
\item[(1)] Ramsey theorem for additive coloring...
\item[(2)] The monadic theory of generalized sums...
\item[(3)] Simple application for decidability...
\item[(4)] The monadic theory of well orderings...
\item[(5)] From orders to uniform-orders...
\item[(6)] Applications of Section 5 to dense orders...
\item[(7)] Undecidability of the monadic theory of the real order...
\end{enumerate}

Let us review the history. Ehrenfeucht \cite{Eh60} proved the decidability
of the first-order theory of order. Gurevich \cite{Gu64} deduced from it the case 
of linear order with one-place predicates. B\"uchi \cite{Bu60} and Elgot \cite{El61}
proved the decidability of the weak monadic theory (i.e., we can quantify over 
finite sets) of (the order of) $\omega$, using automaton theory. B\"uchi continued 
in this direction, in \cite{Bu62}, showing that also the monadic theory (i.e.,
quantification is possible over arbitrary sets) of $\omega$ is decidable; 
and in \cite{Bu65a} he showed the decidability of the
weak
monadic theory of ordinals. 
In   \cite[p. 96]{BS73}% 2021-06-30 10:14 \cite[96]{Bu73} 
he proved the decidability 
of the monadic theory of countable ordinals. 
Rabin \cite{Ra69} proved a very strong and difficult result, 
implying the decidability of 
the monadic theory of countable orders. 
B\"uchi   \cite{BS73} % 2021-06-30 10:20 \cite{Bu73} 
showed the decidability of the monadic 
theory of $\omega_1$ and of $\{\alpha:\alpha<\omega_2\}$. 

Meanwhile La\"uchli \cite{Lc68}, using methods of Ehrenfeucht \cite{Eh59} and 
Fraisse \cite{Fr65} and continuing works of Galvin (unpublished) and La\"uchli 
and Leonard \cite{LcLe66}, proved the decidability of the weak monadic theory of 
order. He did not use automaton theory. Pinus \cite{Pi72} strengthened, somewhat,
those results. Our results have been announced in 
    \cite{Sh:E99}, \cite{Sh:E100} % 2021-06-30 10:17 \cite{Sh73},\cite{Sh73a}.

By our notation La\"uchli used $Th^n_{\bar k}$ only for $\bar k=
\langle 1,1,1,\ldots\rangle$ (changed for the quantification over 
finite sets). 

\par \noindent
Remark: We are not interested here in results without the axiom of 
choice. See Siefkes \cite{Si70} which shows that the result on $\omega$
is provable in ZF. This holds also for $\alpha<\omega^*$. Litman \cite{Li72} % 2021-06-14 02:33 ?}
pointed out some mistakes in \cite[6]{BS73}  % 2021-06-30 10:20 {Bu73} 
(theorems without AC); proved connected 
results, and showed in ZF that $\omega_1$ is always characterizable by a sentence. 

In Section 7 we prove (CH) the undecidability of the monadic theory of 
the real order and of the class of orders, and related problems. It can be 
read independently, and has a discussion on those problems. 
Gurevich finds that our proof works also for the lattice of subsets of a 
Cantor discontinuum, with the closure operation, and similar spaces. 
Hence Grzegorczy's \cite{Gr51} question is answered 
(under CH)\footnote{Gurevich meanwhile has proved more and has a paper in 
preparation.}.

Our work continues \cite{Lc68}, but for well ordering we use 
ideas of B\"uchi and Rabin. We reduce here the decision problem of 
the monadic theories of some (classes of) orders [e.g.,
well orderings; the orders which do not embed $\omega_1$
not $\omega^*_1$] to problems more combinatorial in nature.
So we get a direct proof for the decidability of countable 
orders (answering a question of B\"uchi 
    \cite[p.35]{BS73}  % 2021-06-30 10:21 \cite{Bu73}(p.35)).
Our proof works for a wider class, thus showing that the countable 
orders cannot be characterized in monadic theory, thus 
answering a question of Rabin \cite{Ra69}(p.12). Moreover, there 
are uncountable orders which have the same monadic theory as 
the rationals (e.g., dense Specker order; see \cite{Je03} for % 2021-06-14 15:22 
their existence; and also some uncountable subsets of the reals). 
We also show that the monadic theory of $\{\alpha:\alpha<\lambda^+\}$
is recursive in that of $\lambda$, generalizing results of B\"uchi for 
$\omega$ and $\omega_1$. Unfortunately, even the 
monadic theory of $\omega_2$ contains a statement independent    
of ZFC. For a set $A$ of ordinals, let $F(A)=\{\alpha:\alpha$
is a limit ordinal of cofinality $>\omega,\alpha<\sup A$, 
and $\alpha\cap A$ is a stationary subset of $\alpha\}$. 

Now Jensen \cite{Jn} proved the following: % 2021-07-08 09:10  2021-06-14 02:48 {J1}
\medskip

\begin{theorem}
\label{0.1}
$(V=L)$. A regular cardinal $\kappa$ is weakly compact if and
only if for every stationary $A\subseteqq \kappa$, \st\ $(\forall \alpha\in A)
[\cf (\alpha)=\omega], F(A) \neq \varnothing$.   % 2021-06-14 06:33 cf ()

As the second part is expressible in the monadic theory of order, 
the Hanf number of the monadic theory of order is high. 
Clearly also the monadic theory of the ordinals depends on 
an axiom of large cardinals. 

Now, Baumgartner \cite{B2}  % 2021-06-14 02:52 {Ba1}  a new class of order type
shows that if ZFC+ (there is a weakly 
compact cardinal) is consistent, {\em then} it is consistent with 
ZFC that 
\begin{enumerate}
\item[(*)] for any stationary $A\subseteqq \omega_2$, if 
$(\forall \alpha\in A) [\cf (\alpha)=\omega]$, then $F(A)\neq % 2021-06-14 06:33 cf ()
\varnothing$ (and in fact is stationary). 

So ZFC does not determine the monadic theory of $\omega_2$. 
This partially answers \cite{Bu65}(pp.34-43; p.38, problem 2).

We can still hope that the number of possible such theories is small,
and each decidable, but this seems unlikely. We can also hope to find 
the sentences true in every model of ZFC. A more hopeful project is to 
find a decision procedure assuming $V=L$. We show that for this it suffices 
to prove only the following fact. Let $D_{\omega_2}$ be the filter of closed 
unbounded subsets of $\omega_2$. (Magidor disproves (**) in $V=L$, but it may
still be consistent with ZFC.) 
\item[(**)] if $A\subseteqq \{\alpha<\omega_2:\cf  % 2021-06-14 06:32 cf \ +()
(\alpha)=\omega\}, 
F(A)=B\cup C, A$ is stationary, {\em then} there are $A_1,A_2$, \st\ 
$A=A_1\cup A_2,A_1\cap A_2=\varnothing, A_1,A_2$ are stationary and 
$F(A_1)=B({\rm mod} D_{\omega_2}), F(A_2)=C({\rm mod} D_{\omega_2}).$

We prove, in fact, more: that the monadic theory of $\omega_2$ and the 
first order theory of $\langle \underline{P} (\omega_2)/ D_{\omega_2},\cap,
\cup, F\rangle$ are recursive one in the other.
\end{enumerate}
\end{theorem}

\begin{conjecture}
\label{0.1a}
$(V=L)$. The monadic theory of $\omega_2$ (and even $\omega_n$) is 
decidable. 
\end{conjecture}

\begin{conjecture}
\label{0.2}
($V=L+$ there is no weakly compact cardinal). The monadic theory of 
well orders is decidable. 

La\"uchli and Leonard \cite{LcLe66} define a family $\underline{M}$ of orders 
as follows: It is the closure of $\{1\}$ by 
\begin{enumerate}
\item $M+N$,
\item $M\cdot \omega$ and $M\cdot \omega^*$,
\item $\sum^*_{i<n} M_i$ which is $\sum_{a\in \rationals % 2021-06-14 02:53 Q
    } M_a$ 
and $\{a\in \rationals % 2021-06-14 02:53 Q
    :M_a=M_i\}$ is a dense subset of the rationals, 
and each $M_a\in \{M_i:i<n\}$. 
\end{enumerate}

(See Rosenstein \cite{Ro69} and Rubin \cite{Ru74} for generalization.) % 2021-06-30 10:38 Rt69

L\"auchli \cite{Lc68} proved that every sentence from the weak monadic language 
of order has a countable model if and only if it has a model in $\underline{M}$. 
Easy checking of Section 4 shows this holds also for the monadic language. 
On the other hand, looking at the definition of $\underline{M}$, we can easily see 
that for every $M\in \underline{M}$ there is a monadic sentence $\psi$ \st\ 
$M\models \psi$, and 
$\|N\| \leqq \aleph_0, N\models \psi$ imply $N\cong M$. 

In this way we have a direct characterization of $\underline{M}$. 
\end{conjecture}

\begin{theorem}
\label{0.2a} 
$M\in \underline{M}$ if and only if $M$ is countable and satisfies some 
monadic sentence which is $(\leqq \aleph_0$)-categorical. 

Also for other classes whose decidability we prove, we can find subclasses 
analogous to $\underline{M}$. This theorem raises the following question: 
\end{theorem}

\begin{conjecture}
\label{0.3} 
For every $N\in \underline{M}$ there is a monadic sentence $\psi$ \st\ 
$M\models \psi$ implies that $M$ and $N$ have the same monadic theory. 
(It suffices to prove this for the rational order.)

Related questions are: 
\end{conjecture}

\begin{conjecture}
\label{0.4}
There is a monadic sentence $\psi$ \st\ $\real \models \psi$ and 
$M\models \psi$ imply that $M$ and $R$ have the same monadic 
theory.\footnote{Confirmed by Gurevich}
\end{conjecture}

\begin{conjecture}
\label{0.5} 
There is an order $M$ which has the same monadic theory as $R$, 
but is not isomorphic to $R$.\footnote{Refuted by Gurevich}
\end{conjecture}

\begin{conjecture}
\label{0.6}
There are orders with the same monadic theories, whose completions 
do not have the same monadic theories.\footnote{Confirmed by Gurevich}

The characterization of $\underline{M}$ gives us also
\end{conjecture}

\begin{conclusion}
\label{0.3a}
The question whether a sentence in the first-order (or even 
monadic) theory of order is $(\leqq \aleph_0)$-categorical 
(or $\aleph_0$-categorical) is decidable.

A natural question is whether the monadic theory of ${\mathfrak M}$
is more ``complex'' than that of the ordinals (the orders in ${\mathfrak M}$
are countable unions of scattered types; see Laver \cite[\S3]{Lv71}, 
which includes 
results of Galvin). To answer this, we have the
\end{conclusion}

\begin{definition}
For a model $M$ with relations only, let $M^\sharp$ be the following model:  % 2021-06-14 08:07 ? ==> #
\begin{enumerate}
\item[(i)] its universe is the set of finite sequences of elements of $M$; 
\item[(i)] its relations are 
\begin{enumerate}
\item[(a)] $<$, where $\bar a<\bar b$ means $\bar a$ is a initial segment of 
$\bar b$,
\item[(b)] for each $n$-place predicate $R$ from the language of $M$,
$R^{M^\sharp}=\{\langle\langle a_1,\ldots,a_{m-1},b^1\rangle,  % 2021-06-14 08:08 ? ==> #
\langle a_1,\ldots,a_{m-1},b^2\rangle, \ldots,
\langle a_1,\ldots,a_{m-1},b^n\rangle\rangle:a_i,b^i$ 
are elements of $M$, and $M\models R[b^1,\ldots,b^n]\}$. 
\end{enumerate}

The author suggested a generalization of Rabin's automaton from 
\cite{Ra69}, proved the easy parts: the lemmas on union and intersection, 
and solved the emptiness problem. Then J.Stup elaborated those proofs, 
and proved the complementation lemma. Thus a generalization of the theorem and 
proof of \cite{Ra69} gives
\end{enumerate}
\end{definition}

\begin{theorem}
\label{0.4a}
The monadic theory of $M^\sharp$ is recursive in the monadic   % 2021-06-14 07:16 xob naseq
theory of $M$. 

Thus, using   \cite[\S3]{Lv71} % 2021-06-14 02:55 \cite{Lv1,3} 
    notation, we get, e.g.,
\end{theorem}

\begin{conclusion}
\label{0.5a}
The monadic theory of $\{M:M\in {\mathfrak M},\|M\|
\leqq \lambda\}$ is recursive in the monadic theory of $\lambda$. 

Because by Section 2 the monadic theory of $\sigma_{\lambda^+,\lambda^+}$ 
is recursive in the monadic theory of $\lambda$, by \ref{0.4} the monadic 
theory of $\eta_{\lambda^+,\lambda^+}$ is recursive in the monadic 
theory of $\lambda$, and so we finish, as by \cite[3.2(iv),3.4]{Lv71}  % 2021-06-14 02:58 
$\eta_{\lambda^+,\lambda^+}$ is a universal member of $\{M\in 
{\mathfrak M}:\|M\|\leqq\lambda\}$. 

Also useful are the following (Le Tourneau \cite{To68} % 2021-06-30 10:49 {To71}
proved parts 
(1),(2) at least):\footnote{Le Tourneau only claimed the result. Lately 
also Routenberg and Vinner proved this theorem.}
\end{conclusion}

\begin{theorem}
\label{0.6a}
Let $L$ be a language with one one-place function symbol,  % 2021-06-14 03:02 
    equality 
and one place predicates.
\begin{enumerate}
\item The monadic theory of $L$ is decidable.
\item If a monadic sentence $\psi$ of $L$ has a model, 
it has a model of cardinality $\leqq \aleph_0$.
\item In (2) we can find $n=n(\psi)<\aleph_0$ and a model 
$M$ \st\ $|\{b\in |M|:f(b)=a\}|\leqq n$ for any 
$a\in |M|$. 
\end{enumerate}

This is because, if $M_\lambda$ is the model whose universe is $\lambda$, and 
whose language contains equality only, in $M^\sharp_\lambda$   % 2021-06-14 08:07 # hayah ?
we can interpret a universal 
$L$-model (see Rabin \cite{Ra69}). This implies (1). Note that all 
$M_\lambda$ ($\lambda$ an infinite cardinal) have the same monadic theory. 
This proves (2). For (3) note that if $M_{\aleph_0}\models \psi$, then for 
all big enough $n,M_n\models \psi$.
\end{theorem}

\par \noindent 
Remark (1): \ 
Rabin \cite{Ra69} prove the decidability of the countable Boolean 
algebras, in first-order logic expanded by quantification over ideals. 
By the Stone representation theorem, each countable Boolean algebra 
can be represented as the Boolean algebra generated by the intervals of a 
countable order. By the method of Section 3 we can prove that the theory of 
countable linear orders in monadic logic expanded by quantification 
over such ideals, is decidable, thus reproving Rabin's result. 
(The only points is that methods of Section 2 apply.)
\medskip

\begin{conjecture}
\label{0.7}
The monadic theory of orders of cardinality $\leqq\aleph_1$ 
is decidable when $\aleph_1<2^{\aleph_0}$.
\end{conjecture}

\begin{conjecture}
\label{0.8}
The theory of Boolean algebras of cardinality $<\lambda$ or in first-order 
logic expanded by allowing quantification over ideals is decidable when 
$\lambda\leqq 2^{\aleph_0} (\lambda=\aleph_2\leqq 2^{\aleph_0})$.
\end{conjecture}

\par \noindent
Remark:\ We can prove Conclusion \ref{0.5} by amalgamating the methods 
of Section 4,5, and 6. 
\medskip

\section{Ramsey theorem for additive coloring}

A {\em coloring} of a set $I$ is a function $f$ from the set of 
unordered pairs of distinct elements of $I$, into a finite set $T$
of colors. We write $f(x,y)$ instead of $f(\{x,y\})$, assuming  % 2021-06-14 03:04 . ==>,
usually that $x<y$. The coloring $f$ is additive if for 
$x_i<y_i <z_i\in I$ $ (i=1,2)$.  % 2021-06-14 03:04 rewax=  space
$$
f(x_1,y_1)=f(x_2,y_2); f(y_1,z_1)=f(y_2,z_2)
$$
imply $f(x_1,z_1)=f(x_2,z_2)$. In this case a (partial) 
operation + is defined on $T$, \st\ for $x<y<z\in I,f(x,z)
=f(x,y)+f(y,z)$. 
A set $J\subseteqq I$ is homogeneous (for $f$) 
if there is a $t_0\in T$ \st\ for every $x<y\in J,
f(x,y)=t_0$. 

Ramsey's theorem \cite{Rm29} states, in particular, that if we color 
an infinite set with a finite set of colors, then there is an infinite 
homogeneous subset. This theorem has many generalization and applications. 
It was used in \cite{Bu62} for a coloring which was, in fact, additive. 
Using an idea of Rabin, B\"uchi \cite[12, p.58]{BS73} 
        % 2021-06-30 10:25 {Bu73}  % 2021-06-14 03:05 12,p.58) 
offered an alternative 
proof (using, in fact, additivity) and in \cite[6.2, p.111]{BS73}
    % 2021-06-30 10:27 {Bu73} % 2021-06-14 03:06 (6.2,p.111)
straightforwardly 
generalized it to $\omega_1$ (the result for $\omega_1$ is not true 
for coloring in general). We give the natural extension to arbitrary ordinals
(which is immediate, and included for completeness) and a parallel 
theorem for dense orders. 

\begin{theorem}
\label{1.1}
If $\delta$ is a limit ordinal, $f$ an additive coloring of $\delta$ 
(by a set $T$ of $n$ colors), then there is an unbounded homogeneous 
subset $J$ of $\delta$.
\end{theorem}

\par \noindent 
Remarks:\ \begin{enumerate}
\item If the cofinality of $\delta$ is $\geqq \omega_1$ we can assume 
that if $a, b<c',f(a,c')=f(b,c')$, then $a,b<c\in J$ implies $f(a,c)=f(b,c)$.
\item Instead of $|T|<\aleph_0$, we need assume only $|T|<\cf(\delta)$.
\end{enumerate}
\medskip

\begin{conclusion}
\label{1.2}
Under the condition of \ref{1.1}, there are a closed unbounded subset 
$J$ of $\delta$, and $J_k,J^\ell,1\leqq k,\ell\leqq |T|$ and $t^\ell_k
\in T$ \st\ $J=\cup_k J_k= \cup_\ell J^\ell$, the $J_k$'s 
are disjoint, the $J^\ell$'s are disjoint, and if $a<b\in J,a\in J_k,
b\in J^\ell$ then $f(a,b)=t^\ell_k$.
\end{conclusion}

\begin{theorem}
\label{1.3}
If $f$ is an additive coloring of a dense set $I$, by a finite set $T$ of 
$n$ colors, then there is an interval of $I$ which has a dense homogeneous
subset.
\end{theorem}

\begin{conclusion}
\label{1.4}
Under the hypothesis of \ref{1.3}, there is an interval $(a,b)$ of $I$,
and $(a,b)=\cup^{|T|}_{k=1} J_k=\cup^{|T|}_{\ell=1} J^\ell$ and colors 
$t^\ell_k\in T$ \st\ for $x<y,x\in J_k,y\in J^\ell,f(x,y)=t^\ell_k$. 
\end{conclusion}

\par \noindent 
Remark: We can choose the $J_0,J_k,J^\ell$'s so that they are definable by 
first-order formulas with parameters in the structure $(\delta,<,f)$ (or
$(I,<,f)$). 
\medskip

\par \noindent 
Proof of Theorem \ref{1.1}:\ 
Define: For $x,y\in \delta,x \sim y$ if there is a $z$ \st\ 
$x,y<z<\delta$, and $f(x,z)=f(y,z)$; clearly this implies by the additivity 
of $f$ that for any $z',z<z'<\delta,f(x,z') 
    = f(y,z') $.   % 2021-06-14 03:08 
    It is easy to verify that 
$\sim$ is an equivalence relation with $\leqq |T|$ equivalence classes. 
So there is at least one equivalence class $I$, which is an unbounded subset of 
$\delta$. Let $x_0$ be the first element of $I$. Let, for $t\in T,I_t=
\{y:x_0\neq y\in I,f(x_0,y)=t\}$. 
    Clearly $I-  % 2021-06-14 03:09 -
        \{x_0\}=\cup_{t\in T} I_t$, 
hence for some $s,I_s$ is an unbounded subset of $\delta$. Let 
$\langle a_i:i<\cf (\delta)\rangle$  % 2021-06-14 03:10 ()
    be an increasing unbounded \sq\ of elements of 
$\delta$. Define by induction on $i$ elements $y_i\in I$. 
If for all $j<i (i<\cf (\delta), y_j$ have been defined, let  % 2021-06-14 03:11 (
$y_i<\delta$ be \st\ $y_i>y_j,y_i>a_j,y_i>x_0$ and $f(x_0,y_i)=f(y_j,y_i)$
for any $j<i$, and $y_i\in I_s$. Now $J=\{y_i:i<
    \cf (\delta)\}$ is the desired   % 2021-06-14 03:12 ()
set. Clearly it is unbounded. If $y_j<y_i$ (hence $j<i$) then 
$$
f(y_j,y_i)=f(x_0,y_i)=s.
$$
So $J$ is homogeneous. 
\medskip

\par \noindent 
Proof of Conclusion \ref{1.2}:\ 
If the cofinality of $\delta$ is $\aleph_0$, then the $J$ from \ref{1.1}
is also closed (trivially). So assume 
    $\cf ( \delta)>\aleph_0$, let  % 2021-06-14 03:12 ()
$T=\{t_1,\ldots,t_n\}$, and let $J$, $y_j$ be as defined in the proof of 
\ref{1.1}; and let $J^*$ be the closure of $\{y_{j+1}:j<
    \cf (\delta)\}$.  % 2021-06-14 03:13 
Then $J^*=\{y^j:j<
    \cf (\delta)\}$ is increasing, continuous, and $y^{j+1}= % 2021-06-14 03:13 
y_{j+1}$. Let $J'=\{y^j:j$ is a limit ordinal$\}$, 
$$
\begin{array}{ll}
J_k=\{y^j:j \mbox{ is a limit ordinal, } f(y^j,y^{j+1})=t_k\}, \\ 
J^\ell=\{y^j:j \mbox{ is a limit ordinal, and }
(\forall i<j)  % 2021-06-14 03:15 _ 
    (\exists \alpha) (i<\alpha<j\wedge 
f(y^{\alpha+1},y^j)=t_\ell) \\
\mbox{ but this does not fold for any } \ell'<\ell\}.
\end{array}
$$
Now clearly $J'=\cup_k J_k=\cup_\ell J^\ell$, and if 
$x\in J_k,z\in J^\ell, x<z$ then $x=y^i,z=y^j,i<j,
i,j$ are limit ordinals and there is an $\alpha,i<\alpha<j$, 
\st\ $f(y^{\alpha+1},y^j)=t_\ell$. Hence
$$
\begin{array}{ll}
f(x,z)=f(y^i,y^j)=f(y^i,y^{i+1})+f(y^{i+1},y^{\alpha+1})+
f(y^{\alpha+1},y^j) \\ 
=t_k+f(y_{i+1},y_{\alpha+1})+ t_\ell=t_k+s+t_\ell 
\stackrel{\rm def}{=} t^\ell_k.
\end{array}
$$
Clearly all the demands are satisfied. 
\medskip

\par \noindent 
Proof of Theorem \ref{1.3}: 
Remember that $J\subseteqq I$ is dense in an interval $(a,b)$
if for every $x,y\in I,a<x<y<b$, there is a $z\in J$ \st\ 
$x<z<y$. It is easy to see that if $J\subseteqq I$ is dense 
in an interval $(a,b)$ and $J=\cup^m_{k=  % 2021-06-14 03:17 -
    1} J_k$ $ (m>1)$ then there are   % 2021-06-14 03:17 rewax= space
$k$ and $a',b'$ \st\ $a<a'<b'<b,1\leqq k\leqq m$ and $J_k$ is dense 
in $(a',b')$. 

Define for any $a\in I,J  % 2021-06-14 03:18 I
    \subseteqq I$
$$ F(a,J)=\{t:t\in T, (\forall x> a) (\exists y\in J)
(a<y<x\wedge f (a,y)=t)\}.
$$
Notice, that since $T$ is finite, for any $a\in I$, and any 
$J\subseteqq I$ there is a $b,a<b\in I$ \st:
$$
t\in F(a,J) \mbox{ if and only if there is a } y\in J,a<y<b,f
(a,y)=t.
$$

We define by induction on $m\leqq n 2^n+2$ intervals 
$(a_m,b_m)$, sets $J_m$ dense in $(a_m,b_m)$, and 
(for $m>0$) sets $D_m\subseteqq T$. 

For $m=0$, let $(a_0,b_0)$ be any interval of $I$, and 
$J_0=\{x\in I:a_0<x<b_0\}$. Suppose 
$(a_m,b_m), J_m$ are defined. For any $D\subseteqq T$ 
let $J_m (D)=\{a\in J_m:F(a,J_m)= D\}$. Clearly 
$J_m=\cup_{D\subseteqq T} J_m (D)$ 
and as there are only finitely many 
possible $D$'s $(\leqq 2^n)$, there is an interval 
$(a_{m+1},b_{m+1})$ and $D_{m+1}\subseteqq T$ \st\ 
$J_m(D_{m+1})$ is dense in $(a_{m+1},b_{m+1})$, and 
$a_m<a_{m+1}<b_{m+1}<  % 2021-06-14 03:18 <
    b_m$. Let $J_{m+1}=
(a_{m+1},b_{m+1})\cap J_m (D_{m+1})$. Clearly 
$J_m\supseteqq J_{m+1}$, and $m>k$ implies 
$J_k\supseteqq J_m$, and $(a_m,b_m)$ 
is a subinterval of $(a_k,b_k)$.

As there are only $\leqq 2^n$ possible
$D_m$, there are a $D\subseteqq T$ and $0\leqq m_0
<\ldots< m_n \leqq  n 2^n+1$ \st\ $D_{m_i+1}=D$. 
Define, for $0\leqq k\leqq n,a^k=a_{m_k},
b^k=b_{m_k}, J^k=J_{m_k}$.\footnote{In fact 
$D_m(T)\supseteqq D_m (T)$, hence we can replace 
$n 2^n+2$ by $n^2+2$.}

It is easy to check that if $0\leqq k<l\leqq n,x\in 
J^\ell$ then $x\in J_{m_\ell}\subseteqq 
J_{m_{k+1}}$, hence $F(x,J_k)=F(x,J_{m_k})=
D_{m_{k+1}}=D$. It is clear that $J^0\supseteqq J^1
\supseteqq\ldots \supseteqq J^n$.

Choose $x_0\in J^n$. Then there is $x_1,x_0<x_0<x_1<b^n$, 
\st\ $x_0<y<x_1,y\in J^0$ implies $f(x_0,y)\in F(x_0,J^0)=D$. 
Hence $t\in D$ if and only if there is $y\in J^{n-1},x_0<
y<x_1,f(x_0,y)=t$, if and only if there is $y\in J_0,
x_0<y<x_1,f(x_0,y)=t$. Clearly     
$$
J^n\cap (x_0,x_1)=\cup_{t\in T} \{y:y\in J^n,
x_0<y<x_1,f(x_0,y)=t\}.
$$
Hence there are $a,b,t_0$ \st\ $x_0<a<b<x_1$ and
$$
J^*=\{y:y\in J^n,a<y<b,f(x_0,y)=t_0\}
$$
is dense in $(a,b)$. Clearly $t_0\in D$.  

It is easy to check that for $t,s\in D,t+s$
is defined and $\in D$, so for $t\in D, m\geqq 1$ 
defined $mt \in T$, by induction on $m:1t=t, 
(m+1) t=mt+t$. As $T$ has $n$ elements, 
$1t_0,2t_0,\ldots, (n+1) t_0$ cannot be pairwise 
distinct. So there are 
$i,j,1\leqq i < (i+j) \leqq n+1$ \st\ $it_0
=(i+j) t_0$. Define
$$
J=\{y:a<y<b, f(x_0,y)=j t_0,y\in J^{n-j+1}\}.
$$
We shall show that $J$ is the desired set.
\begin{enumerate}
\item[(I)] $J$ is dense in $(a,b)$.

Suppose $a<a'<b'<b$, and we shall find $z\in J,a'
<z<b'$. As $J^*$ is dense in $(a,b)$ there are $z^n\in 
J^*\subseteqq J^n,a'<z^n<b'$. We define by downward 
induction $z^k$ for $n\minus   j+1\leqq k\leqq n$ \st\ 
$z^k\in J^k,a'<z^k<b'$. For $k=n,z^k$ is defined. 
Suppose $z^{k+1}$ is defined, then as 
$z^{k+1}\in J^{k+1}$ is follows that 
$F(z^{k+1},j^k)=D$. As $t_0\in D$ there is 
$z^k\in J^k$, \st\ $z^{k+1}<z^k<b'$ 
and $f(z^{k+1},z^k)=t_0$. Clearly
$$
x<z^n<z^{n-1}<\ldots<z^{n-j+1},
$$
$$
f(x_0,z^n)=t_0,\quad f(z^{i+1},z^i)=t_0.
$$
Hence $f(x_0,z^{n-j+1})=t_0+\ldots
+t_0=j t_0$, so $z^{n-j+1}\in J,a'<
z^{n-j+1}<b'$.
\item[(II)] $J$ is homogeneous. 

Suppose $a<y<z<b,y,z\in J$. Then $y\in J^{n-j+1}$. 
Now define by downward induction $y^k\in J^k$ for 
$0\leqq k\leqq i,y\leqq y^k<z$. Let $y^i=y(y^i\in J^i$ 
because $y^i=y\in J^{n-j+1}$, and as $i+j\leqq n+1,i
\leqq n\minus   j+1$ hence $J^{n-j+1}\subseteqq J^i).$ If 
$y^{k+1}$ is defined then $F(y^{k+1},J^k)=D$, hence 
there are $y^k\in J^k, y^{k+1}<y^k<z$ \st\ 
$f(y^{k+1},y^k)=t_0$. It follows that $x_0<
y=y^i<y^{i-1}<\ldots< y^0<z$ and  
$$
f(y^k,y^{k-1})=t_0.
$$
Hence 
$$
f(y,y^0)=f(y^i,y^0)=it_0.
$$
So 
$$
f(y,z)=f(y,y^0)+f(y^0,z)=it_0+f(y^0,z)
$$
$$
\quad \quad =(i+j) t_0+f(y^0,z)=jt_0+it_0+f(y^0,z)
$$
$$
\quad \quad =f(x_0,y)+f(y,y^0)+f(y^0,z)=f(x_0,z)=jt_0.
$$
This proves the homogeneity of $J$. 
\end{enumerate}
\medskip

\par \noindent
Proof of Conclusion \ref{1.4}:\ 
Let $(a,b),J$ and $t_0$ be as in the proof of \ref{1.3}.
Let $T=\{t_1,\ldots,t_n\}$. Let 
$$
J_k=\{y:y\in (a,b),t_k\in F(y,J),t_1,\ldots, t_{k-1}
\notin F(y,J)\},
$$
$$
J^\ell=\{y:y\in (a,b),t_\ell\in F' (y,J),t_1,\ldots, 
t_{\ell-1}\notin F' (y,J)\}
$$
where $F'$ is defined just as $F$ is, but for the 
reversed order. 

Clearly $(a,b)=\cup_k J_k=\cup_\ell J^\ell$. Suppose 
$x<y,x\in J_k,y\in J_\ell$. Then we can find 
$x',y' $  $x<x'<y'\in J$, \st\ $f(x,x')=t_k,f(y',y)=t_\ell$.   % 2021-06-14 03:20 rewax
Hence 
$$
f(x,y)=f(x,x')+f(x',y')+f(y',y)=t_k+t_0+t_\ell \stackrel{{\rm def}}{=}
t^\ell_k.
$$
\medskip

\section{The monadic theory of generalized sums}

Feferman and Vaught \cite{FeVa59}  % 2021-06-14 03:20 {FV1} 
    proved that the first order theory 
of sum, product, and even generalized products of models depends 
only on the first-order theories of the models. 
Their theorem has generalizations to even more general products 
(see Olmann)   % 2021-06-14 03:30 xob xob % 2021-06-14 03:30 \cite{Ol,1}) 
and to suitable infinitary 
languages ($L_\alpha$, see Malitz \cite{Mal71}).  % 2021-06-14 03:29   % 2021-07-19 14:26 Ma71  % 2021-07-20 09:24 Ma71==>Mal71

        % 2021-07-20 09:29 Ml71==> Mal71
On the other hand, it is well-known 
that for second order theory this is false even for sum 
(as there is a sentence true in the sum of two models if and only 
if they are isomorphic, for fixed finite language, of course). Also 
for monadic (second-order) theory this is false for products of 
models (there is a sentence true in a direct product of two 
models of the theory of linear order if and only if the 
orders are isomorphic). We notice here that the monadic 
theory of generalized sum depends only on the monadic theories 
of the summands and notice also generalization of known refinement 
(see Fraiss\'e \cite{Fr65}). We can prove them using natural generalization 
of Ehrenfeucht games (see \cite{Eh60}). L\"auchli \cite{Lc68} uses some 
particular cases of those theorems for the weak monadic theory. As 
there is no new point in the proofs, we skip them. We should notice 
only that a subset of sum of models is the union of subsets of 
the summands. The results of \cite{FeVa59} can be applied directly by 
replacing $M$ by $(|M|\cup \underline{P}(M),M,\in)$.

\begin{notation}  
$L$ will be first-order language with a finite number of symbols,
$L^M$ the corresponding monadic language, $L(M)$ the first-order, 
language corresponding to the model $M$, the universe of $M$, is 
$|M|$. Let $x,y,z$ be individual variables; $X,Y,Z$ set variables;
$a,b,c$ elements; $P,Q$ sets; 
$\underline{P}(M)=\{P:P\subseteqq |M|\}$. Bar denotes that this is a finite sequence, 
e.g., $\bar a$; $\ell(\bar a)$ its length, 
$\bar a=\langle \ldots,a_i,\ldots\rangle_{i<\ell(\bar a)}$, and let 
$\bar a(i)=a_i$. We write $\bar a\in A$ instead of $a_i\in A$ and 
$\bar a\in M$ instead of $\bar a\in |M|$. $K$ is a class of 
$L(K)$ models $(L(K)=L(M)$ for any $M\in K$). Let 
$$
K^m=\{(M,\bar P):\bar P\in \underline{P}(M)^m\},K^\infty
=\cup_{m<\omega} K^m.
$$
Let $k,\ell,m,n,p,q,r$ denote natural numbers. 
\end{notation}

\begin{definition}
\label{2.1}
For any $L$-model $M, \bar P\in \underline{P} (M), 
\bar a\in |M|,\Phi$ a finite set of formulas $\varphi
(X_1,\ldots,x_1,\ldots)\in L$, a natural number $n$, and 
a \sq\ of natural numbers $\bar k$ of length 
$\geqq n$, define 
$$
t=th^n_{\bar k} ((M,\bar P,\bar a),\Phi)
$$
by induction on $n$: 

For $n=0$:
$$
t=\{\varphi(X_{\ell_1},\ldots,x_{j_1},\ldots):
\varphi(X_1,\ldots,x_1,\ldots)\in \Phi,M\models \varphi
[P_{\ell_1},\ldots,a_{j_1},\ldots]\}.
$$
For $n=m+1$:
$$
t=\{th^m_{\bar k}(M,\bar P,\bar a\conc \bar b):
\bar b\in |M|^{\bar k(m)}\}.
$$
\end{definition}

\begin{definition}
\label{2.2}
For any $L$-model $M,\bar P\in \underline{P} (M)$, a finite set 
$\Phi$ of formulas $\varphi(X_1,\ldots,x_1\ldots)\in L,n,\bar k$ 
of length $\geqq n+1$, define $T=Th^n_{\bar k} (M,\bar P),\Phi)$ by 
induction on $n$:

For $n=0$:
$$
T=th^1_{\bar k} ((M,\bar P),\Phi).
$$
For $n=m+1$:
$$
T=\{Th^m_{\bar k} ((M,\bar P\conc \bar Q),\Phi):\bar Q\in 
\underline{P} (M)^{\bar k(n)}\}.
$$
\end{definition}

\par \noindent
\begin{enumerate}
\item[(1)] If $\Phi$ is the set of atomic formulas we shall omit 
it and write $Th^n_{\bar k} (M,\bar P)$.
\item[(2)] We always assume $\bar k(i)\geqq 1$ for any $i<\ell(\bar k)$,
and $\bar k(0)\geqq m_R $ % 2021-06-14 03:38 n$ 
    if $R\in L(M)$ is $m_R$-place.
\item[(3)] If we write $\bar k(i)$ for $i\geqq \ell(\bar k)$, then 
we mean 1, and when we omit $\bar k$ we mean $\langle {\rm max} 
    \{m_R: R \in L(M) \}   % 2021-06-14 03:39 m_R
,1,\ldots\rangle$. 
\item[(4)] We could have mixed Definition \ref{2.1}, and \ref{2.2},
and obtained a similar theorem which would be more refined. 
\end{enumerate}
\medskip

\begin{lemma}
\label{2.1a}
\begin{enumerate}
\item[(A)] For every formula $\psi(\bar X)\in L^M (M)$  there is an $n$
\st\ from $Th^n_{\bar k} (M,\bar P)$ we can find effectively whether 
$M\models \psi[\bar P]$.  
\item[(B)] For every $L,\bar k,n,\Phi\subseteqq L$, and $m$ there is a 
set $\Psi=\{\psi_\ell(\bar X):\ell<\ell_0(<\omega),\ell(\bar X)=m\}
(\psi_\ell\in L^M)$ \st\ for any $L$-models $M,N$ and  % 2021-06-14 07:12 L_M
$\bar P\in \underline{P}(M)^m,\bar Q\in \underline{P} (N)^m$ the following 
hold: 
\begin{enumerate}
\item $Th^n_{\bar k}((N,\bar Q),\Phi)$ can be computed from 
$\{\ell<\ell_0:N\models \psi_\ell [\bar Q]\}$.
\item $Th^n_{\bar k} ((N,\bar Q),\Phi)=Th^n_k((M,
\bar P),\Phi)$ if and only if for any $\ell<\ell_0, 
M\models \psi_\ell [\bar P]\Leftrightarrow N \psi_\ell
[\bar Q]$. 
\end{enumerate}
\end{enumerate}
\end{lemma}

\par \noindent 
Proof: Immediate. 
In (A) it suffices to take for $n$ the quantifier depth 
of $\psi$.
\medskip

\begin{lemma}
\label{2.2a}
\begin{enumerate}
\item[(A)] For given $L,n,m, \bar k$, each $Th^n_{\bar k}(M,\bar P)$ is  % 2021-06-14 07:13 ,
hereditarily finite, and we can compute the set of formally possible  
$Th^n_{\bar k} (M,\bar P),\ell(\bar P)=m,M$ an $L$-model. 
The same holds for $\Phi$. 
\item[(B)] If $\bar \ell(0)\geqq \bar k(0),1=p_0<p_1<p_2<
\ldots<p_n\leqq m$ and for 
$1\leqq i\leqq n , \text{ }\bar k(i)\leqq 
\sum_{p _{i-1}\leqq j\leqq p_i} \bar \ell (j)$ then from   % 2021-06-14 07:11 P==> p
$Th^m_{\bar \ell} ((M,\bar P),\Phi)$ we can effectively compute 
$Th^n_{\bar k} ((M,\bar P),\Phi)$. 
\item[(C)] For every $n,\bar k,\bar \ell$ we can compute $m$ 
\st\ from $Th^m_{\bar \ell} ((M,\bar P),\Phi)$ we can effectively 
compute $Th^n_{\bar k} ((M,\bar P),\Phi)$. 
\item[(D)] Suppose in Definition \ref{2.2} we make the following changes: 
We restrict ourselves to partition $\bar P$, and let 
$\bar Q$ be a partition refining $\bar P$, which divides each 
$P_i$ to $2^{\bar k(m)}$ parts. What we get we call $p
Th^n_{\bar k} ((M,\bar P),\Phi)$. Then from $p Th^n_{\bar k}
((M,\bar P),\Phi)$ we can effectively compute $Th^n_{\bar k} 
((M,\bar P),\Phi)$, and vice versa. 
\item[(E)] Let $K,n,\Phi$ be given. If for every $\bar k$ there is 
an $\bar \ell$ \st\ for every $m,M,N\in K^m$, 
$$
Th^n_{\ell}  % 2021-06-14 07:10 
    (M,\Phi)=Th^n_{\bar \ell} (N,\Phi) \Rightarrow 
Th^{n+1}_{\bar k} (M,\Phi)=Th^{n+1}_{\bar k} (N,\Phi)
$$
then for every $m,\bar k$ there is an $\bar \ell$ \st\ for 
any $n' , % 2021-06-14 07:09 m,
    M,N\in K^m$
$$
Th^n_{\bar \ell} (M,\Phi)=Th^n_\ell (N,\Phi) \Rightarrow 
Th^{n'}  % 2021-06-14 07:08 m
_{\bar k} (N,\Phi)= Th^{n'}_{\bar{ k }}(M, \Phi ).
$$
\end{enumerate}
\end{lemma}

\par \noindent 
Remark: This is parallel to elimination of quantifiers. 

(F)\ In (E), if in the hypothesis $\bar \ell$ can be found effectively 
from $\bar k$ then in the conclusion, $\bar \ell$ can be found effectively 
from $m,\bar k$. If in addition $\{Th^n_{\bar k} (M,\Phi):
M\in K^m\}$ is recursive in $\bar k,m$ then 
$\{Th^p_{\bar k} (M,\Phi):M\in K\}$ is recursive in $p,\bar k$. 
\medskip

\par \noindent 
Proof: Immediate. 

The following generalizes the ordered sum of ordered sets (which will be our
main interest) to the notion of a generalized sum of models. 
(Parts (1),(2),(3) of the definition are technical 
preliminaries.)
\medskip

\begin{definition}
\label{2.3}
Let $L_1,L_2,L_3$ be first-order languages, $M_i$
an $ L_1$-model
    % 2021-06-14 07:06 and $L_i$-
 (for $i\in |N|),N$ an $L_2$-model, and we shall define the $L_3$-model 
$M=\sum^\sigma_{i\in |N|} M_i$ (the generalized sum of the 
$M_i$'s relative to $\sigma$).\footnote{We assume, of course, that the 
$|M_i|$'s are pairwise disjoint.}
\begin{enumerate}
\item An $n$-condition $\tau$ is a triple 
$\langle E,\Phi,\Psi\rangle$ where: 
\begin{enumerate}
\item[(A)] $E$ is an equivalence relation on 
$\{0,1,\ldots,n-1\}$.
\item[(B)] $\Phi$ is a finite set of formulas of the form 
$\varphi(x_{j_1},\ldots,x_{j_k})$ where $j_1,\ldots, 
j_k$ are $E$-equivalent and $<n$; and $\varphi\in L_1$. 
\item[(C)] $\Psi$ is a finite set of formulas of the form 
$\psi(x_{j_1},\ldots,x_{j_k})$ where $j_1,\ldots,j_k<
n,\psi\in L_2$. 
\end{enumerate}
\item If $a_0,\ldots,a_{n-1}\in \bigcup_{i\in |N|} M_i,\tau=
\langle E,\Phi,\Psi\rangle$ is an $n$-condition, $a_\ell\in 
M_{i(\ell)},$ then we say 
$\langle a_0,\ldots,a_{n-1}\rangle$ satisfies $\tau$ if:
\begin{enumerate}
\item[(A)] $i(\ell)=i(m)\Leftrightarrow \ell Em$;
\item[(B)] $\varphi(x_{j_1},\ldots,x_{j_k})\in \Phi
\Rightarrow    % 2021-06-14 07:05 <===
M_{i(j_1)} \models \varphi 
[a_{j_1},\ldots,a_{j_k}]$;
\item[(C)] $\psi (x_{j_1},\ldots,x_{j_k})\in \Psi
\Rightarrow  N \models \psi [i(j_i),\ldots,i(j_k)].$  % 2021-06-14 07:05 <===
\end{enumerate}
\item The rule, $\sigma$ is $\langle L_1,L_2,L_3,\sigma^*\rangle$ 
where $\sigma^*$ is a function whose domain is the set of 
predicates of $L_3$; if $R$ is an $n$-place predicate in $L_3,
\sigma^*(R)$ will be a finite set of $n$-conditions. 
\item $M =  % 2021-06-14 07:03 -
    \sum^\sigma_{i\in|N|} M_i$ is an $L_3$-model,  % 2021-06-14 07:04 s,
        whose universe is 
$\cup_{i\in |N|} |M_i|$, and for every predicate 
$\real \in L_3,R^M=\{\langle a_0,\ldots,a_{n-1}\rangle$ satisfies some 
$\tau\in \sigma^*(R)\}$. 

Let $\Phi(\sigma) \text{ } (\Psi(\sigma))$ be the set of all formulas   % 2021-06-14 07:02 xob
$\varphi_j\in L_1 (\sigma) \text{ }  (\psi_p\in L_2(\sigma))$ appearing in 
the $\sigma(R)$'s, $R\in L_3(\sigma)$, and the equality. 
\end{enumerate}
\end{definition}

\par \noindent 
Remarks: \begin{enumerate}
\item We use the convention that 
$  
\sum^\sigma_{i\in N} (M_i,\bar P^i)= (\sum^\sigma_{i\in N} M_i,
\cup_{i\in N} \bar P^i)
$  % 2021-06-14 03:42 $
where for $\bar P^i=\langle P^i_1,\ldots, P^i_m\rangle, 
\bigcup_i \bar P_i=\langle \bigcup_i P^i_1,\ldots,
\bigcup_i P^i_m\rangle$.
\item We could have defined the sum more generally, by allowing the 
universe and the equality to be defined just as the other relations. 
\end{enumerate}
\medskip

\begin{lemma}
\label{2.3a}
For any $\sigma,n,m,\bar k$, if for $\ell=1,2,
\bar P^\ell_1\in \bar P (M^\ell_i)^m$ and for every $i\in N$,
$$
Th^n_{\bar k} ((M^1_i,\bar P^1_i),\Phi (\sigma))=
Th^n_{\bar k} ((M^2_i,\bar P^2_i,\bar P^2_i),\Phi(\sigma)),
$$
then 
$$
Th^n_{\bar k} (\sum^\sigma_{i\in N} (M^1_i,\bar P^1_i))=
Th^n_{\bar k} (\sum^\sigma_{i\in N} (M^2_i,\bar P^2_i)),
$$
\end{lemma}

\begin{theorem}
\label{2.4}
For any $\sigma,n,m,\bar k$ we can find an $\bar r$ \st: 
if $M=\sum^\sigma_{i\in N} M_i,t_i=Th^n_{\bar k} ((M_i,
\bar P_i),\Phi(\sigma))$, and $Q_t=\{i\in N:t_i=t\},
\ell(\bar P_i)=m$, then from $Th^n_{\bar r} ((N,\ldots,Q_t,
\ldots),\Psi (\sigma))$ we can effectively compete 
$Th^n_{\bar k} (M,\bigcup_i \bar P_i)$ 
(which is uniquely determined).
\end{theorem}

\begin{definition}
\label{2.4a}
\begin{enumerate}
\item  For a class $K$ of models
$$
Th^n_{\bar k} (K,\Phi)=\{Th^n_{\bar k} (M,\Phi): M\in K\}.
$$
\item The monadic theory of $K$ is the set of monadic sentences true in 
every model in $K$. 
\item For any $\bar \sigma,K_1,K_2,$ let 
$C\ell^{\bar \sigma} (K_1,K_2)$ be the minimal class $K$
\st\ 
\begin{enumerate}
\item[(A)] $K_1\subseteqq K$,
\item[(B)] if $j<\ell (\bar \sigma), M_i\in K,N\in K_2$ 
then $\sum^{\bar \sigma(i)}_{i\in |N|} M_i\in K$. 
\end{enumerate}
\end{enumerate}
\end{definition}

\begin{conclusion}
\label{2.5}
Suppose $\bar \sigma,n,\bar k,m$ are given. 
$L_1(\sigma_i)=L_3(\sigma_i)=L,L_2(\sigma_i)=L_2;
L,L_2$ are finite and each $\Psi(\sigma_i),\Psi(\sigma_i)$
is a set of atomic formulas. There is an $\bar r$ \st\ for 
every $K_1,K_2,$ from $Th^n_{\bar r} 
(K^{\bar r(n+1)}_2), Th^n_{\bar k} (K_1^m)$ we can effectively 
compute $Th^n_{\bar k}(K^m)$ where 
$K=C\ell^{\bar \sigma} (K_1,K_2)$ (remember 
$K^m_1=\{(M,\bar P):M\in K_1,\bar P\in 
\underline{P} (M)^m)$ ($K_1$ should be a class of 
$L$-models, $K_2$ a class of $L_2$-models).
\end{conclusion}

\par \noindent 
Proof:\ For every $j<\ell (\bar \sigma)$ let $\bar r^j$ 
relate to $\bar \sigma (j),n,\bar k,m$ just as 
$\bar r$ relates to $\sigma,n,k,m$ in Theorem \ref{2.4}. 
Now choose an $\bar r$ \st\ for every $\ell\leqq n, 
\bar r (\ell) \geqq r^j (\ell)$. 

Let $T$ be the set of formally possible 
$Th^n_{\bar k} (M,\bar P)$, for $M$ and $L$-model,
$\ell(\bar P)=m$, and we can define $r(n+1)=|T|$. 
Let $T=\{t(0),\ldots,t(p-1)\}$ 
(so $p=|T|=r(n+1))$. 

Clearly, by the definition of $\bar r^j$, and by 
(a trivial case of) \ref{2.2}(B), if 
$M=\sum^{\bar \sigma(j)}_{i\in N} M_i,t_i=Th^n_{\bar k}
(M_i,\bar P_i),Q_\ell=\{i\in N:t_i=t(\ell)\},\ell(\bar P_i)=m$,
then from $t=Th^n_{\bar r} (N,\ldots, Q_1,\ldots)_{\ell<p}$ 
we can effectively compute $Th^n_{\bar k} (M,\bigcup_i,\bar P_i)$, 
and denote it by $g(t)$. 

Now define by induction on $\ell,T_\ell\subseteqq T$.

Let $T_0=Th^n_{\bar k} (K_\ell^m)$, and if $T_q$ is defined let 
$T_{q+1}$ be the union of $T_q$ with the set of 
$t\in T$ satisfying the following condition:
\begin{enumerate}
\item[(*)] There is a $t^*\in Th^n_{\bar r} 
(K^{r(n+1)}_2)$ \st\ $t=g(t^*)$, and if $t^*$ implies that 
$Q_\ell$ is not empty, then $t(\ell)\in T_q$.
\end{enumerate}
\medskip

\par \noindent 
Remark:\ Clearly if $t^*=Th^n_{\bar r} (N,\ldots,Q_\ell,\ldots)$ then 
from $t^*$ we can compute $Th^0_{\bar r} (N,\ldots,Q_\ell,\ldots)$ and 
hence know whether $Q_\ell\neq \varnothing$.

Clearly $T_0\subseteqq T_1\subseteqq T_2,\ldots\subseteqq T$ so, as 
$|T|=p$, for some $q\leqq p, T_q=T_{q+1}$.

Now let
$$
K_*=\{M\in K: \mbox{ for every } \bar P\in (\underline{P} 
(|M|)^m Th^n_k (M,\bar P)\in T_q\}.
$$
Clearly $Th^n_{\bar k} (k^m_*)\subseteqq T_q$, and we can 
effectively find $T_q$. Now if $N\in K_2, M_i\in K_*$ for $i\in N$, and 
$M=\sum^{\sigma(j)}_{i\in N} M_i$, then for any 
$\underline{P} \in \bar P (|M|)^m, Th^n_{\bar k} (M,P)\in 
T_{q+1}=T_q$ by the definition of $T_{q+1}$, and $M\in K$ 
by the definition of $K$, hence $M\in K_*$. As clearly 
$K_1\subseteqq K_* \subseteqq K$, by the definition of $K=
C\ell^{\bar \sigma} (K_1,K_2)$ necessarily $K_*=K$. So it 
suffices to prove that $Th^n_{\bar k} 
(K^m_*)\supseteqq T_\ell$. (Take $\ell=q$.) This is done by 
induction on $\ell$. 
\medskip   

\begin{lemma}
\label{2.6}
If $M$ is a finite model, then for any $\Phi,n,\bar k$ we can 
effectively compute $Th^n_{\bar k} (M,\Phi)$ from $M$. 
\end{lemma}

\begin{remark}
\label{2.7}
Naturally we can ask whether we can add to (or replace the) 
monadic quantifiers (by) other quantifiers, without essentially 
changing the conclusions of this section. 
It is easily seen that, e.g., the following quantifiers suitable:
\begin{enumerate}
\item $(\exists^f X)$ --there is a finite set $X$
\item $(\exists^\lambda X)$ --there is a set $X,|X|<\lambda$ ($\lambda$
a regular cardinal). when dealing with ordered sums of linear order, also
\item $(\exists^{wo} X)$ --there is a well-ordered set $X$
\item $(\exists_\lambda X)$ --there is a set $X$, with no increasing not 
decreasing sequence in it of length $\lambda$ ($\lambda$ a regular 
cardinal).
\end{enumerate}

If we add some of those quantifiers, we should, in the definition of 
$Th^0_n ((M,\bar P),\Phi)$ state which Boolean combinations of the $P_\ell$'s 
are in the range of which quantifiers. If we e.g., replace the monadic quantifier 
by $(\exists^\lambda X)$, we should restrict the $P$'s to sets of cardinality 
$<\lambda$. 

Another possible generalization is to generalized products. Let $M=
\prod^\sigma_{i\in N} M_i$ (where $L(M_i)=L_1(\sigma), L(N)=L_2(\sigma),
L(M)=L_3(\sigma))$ 
means: $|M|=\prod_{i\in N} |M_i|$, and if 
$f_1,\ldots, f_n\in M,M\models R [f_1,\ldots,f_n]$ if and only if 
$N\models \psi_R [\ldots, P_\ell,\ldots]$ where 
$$
P_\ell=\{i\in N:M_i\models \varphi^R_\ell [f_1(i),\ldots,
f_n(i)]\}
$$
(and $\varphi_\ell$ is a first order sentence from $L_1
(\sigma),\psi_R$ a monadic sentence from $L_3(\sigma)$). Then,
of course, we use $Th^n_{\bar k} (N,\underline{P}), th^n_{\bar k} 
(M,\bar a)$. All our % 2021-06-14 06:58 out 
    theorems generalize easily, but still no 
application was found. 

If not specified otherwise, we restrict ourselves to the class 
$K_{\rm ord}$ of models of the theory of order (sometimes with 
one-place relations which will be denoted, e.g., 
$(M,\bar P)$). $\sigma=\sigma_{\rm ord}$ is the ordered sum of 
ordered sets and is omitted. Therefore $\Psi(\sigma)$ and   % 2021-06-14 06:58 \psi 
$\Phi(\sigma)$ are the set of atomic formulas. For the sum of 
two orders we write $M_1+M_2$. The ordinals, the reals 
$\real $, and the rationals $\rationals $ have their natural orders.   % 2021-06-14 06:59 R Q
If $M=\sum_{i\in |N|} M_i$ we write $Th^n_{\bar k} 
(M,\bar{ P }    ) % 2021-06-14 07:00 \underline{P})
=\sum_{i\in |N|} Th^n_{\bar k} (M_i,\bar P_i)$ where 
$\bar P=\bigcup_i \bar P_i$. Let $T(n,m,\bar k)$ be the set of formally 
possible $Th^n_{\bar k} (M,\bar P),M$ an order,   % 2021-06-14 07:00 th 
$\ell(\bar P)=m$. 
\end{remark}

\begin{corollary}
\label{2.8}
For any $n,m,\bar k$ there is $\bar r=\bar r (n,m,\bar k)$ \st\ if 
$P_t=\{i\in N:t_i=t\}$ for $t\in T(n,m,\bar k)$ then 
$\sum_{i\in N} t_i$ can be effectively 
computed from $Th^n_{\bar r} (N,\ldots,P_t,\ldots)$.
\end{corollary}

\section{Simple application for decidability}

Using Section 2 we shall prove here some theorems, most of them known. 
We prove the decidability of the theories of the finite orders, the countable
ordinals \cite{BS73}  % 2021-06-30 10:26 {Bu73}, 
and show that from the monadic theory of $\lambda$ we can compute 
effectively the monadic theory of $K=\{\alpha:\alpha<\lambda^+\}$ 
(this was shown for $\lambda=\omega,\lambda=\omega_1$ in \cite{BS73}
        % 2021-06-30 10:26 {Bu73}. 
We do not try to prove the results on definability and elimination of 
quantifiers. For finite orders this can be done and the method becomes similar 
to that of automaton theory. For $\omega,\{\alpha:\alpha<\omega_1\},\omega_1$ this
can be done by using the previous cases (e.g., for $\omega$ using the result on 
the finite orders). We can prove the decidability of the weak monadic theory 
(with $\exists^f $ % 2021-06-14 06:53 t$ 
    only) of the $n$-successors theory by the method of this section 
(Doner \cite{Do65} proved it). It would be very interesting if we could have 
proved in this way that the monadic theory of the 2-successor theory is 
decidable (Rabin \cite{Ra69} proved it). 

In order to use Section 1 we should note 

\begin{lemma}
\label{3.1} 
For any $m,\bar k, (N,\bar P)$, the coloring $f^n_{\bar k}$ on $N$ is additive 
where 
$$
f^n_{\bar k} (a,b)=Th^n_{\bar k} ((N,\bar P)\rest [a,b)),
$$
where
$(N,\bar P)\rest [a,b)$  % 2021-06-14 06:53 (
    is a submodel of $(N,\bar P)$ with the universe 
$[a,b)=\{x\in N:a\leqq x<b\}$.   
\end{lemma}

Proof:\ By lemma \ref{2.3a}.

Let us list % 2021-06-14 06:54 limit 
    some immediate claims. 
\medskip

\begin{lemma}
\label{3.2}
\begin{enumerate}
\item[(A)] If for any $n,\bar k$ we can compute effectively 
$Th^n_{\bar k} (K)$, then the monadic theory of $K$ is 
decidable; and vice-versa.
\item[(B)] If the monadic theory of $K$ is decidable then so is the 
monadic theory of $K'$ where $K'$ is the class of:
\begin{enumerate}
\item[(i)] submodels of $K$,
\item[(ii)] initial segments of orders from $K$,
\item[(iii)] orders which we get by adding (deleting) first 
(last) elements from orders of $K$, 
\item[(iv)] converses of orders from $K$,
\item[(v)] $(M,\bar P),M\in K,\bar P\in \underline{P}
(M)^m$. 
\end{enumerate}
\end{enumerate}
\end{lemma}

Proof: \ Immediate.
\medskip

\begin{theorem}
\label{3.3}
The monadic theory of the class $K_{\rm fin}$ of finite orders is 
decidable.
\end{theorem}

Proof:\ Let $K_n$ be the class of orders of cardinality $n$; up to 
isomorphism $K_n$ has only one element, $n$. Hence by Lemma \ref{2.6} we 
can compute $Th^n_{\bar k} (K_i)$. Hence by Conclusion \ref{2.5}, for every 
$n,\bar k$ we can compute $Th^n_{\bar k} (K)$ where $K= C\ell (K_1,K_2)$. 
But clearly $K$ is the class of finite orders. So by \ref{3.2}(A) we finish.
\medskip

\begin{theorem}
\label{3.4}
The monadic theory of $\omega$ is decidable. 
\end{theorem}

Proof:\ We shall compute $\{Th^n_{\bar k} (\omega,\bar P):
\bar P\in \underline{P} (\omega)^m\}$ by induction on $n$, 
for every $\bar k,m$ simultaneously. 

For $n=0$ is it easy. 

Suppose we have done it for $n -  % 2021-06-14 06:54 - 
    1$ and we shall do it for 
$n,m,\bar k$. By the induction hypothesis we can compute $Th^n_{\bar \ell}
(\omega)$ for every $\bar \ell$, in particular for $\bar r=\bar r
(n,m,\bar k)$ (see \ref{2.8}). Now for any $M=(\omega,P_1,\ldots,
P_m)$, by \ref{1.1} we can find an $f^n_{\bar k}$-homogeneous set 
$\{a_i:i<\omega\} (a_i<a_{i+1})$. So letting 
$$ t=T^n_{\bar k} ((\omega,\bar P)\rest [0,a_0)),\ s=Th^n_k
((\omega,\bar P)\rest [a_i,a_j))\ \mbox{for } i<j;
$$
we have 

\quad $Th^n_{\bar k}(\omega,\bar P)=Th^n_{\bar k} ((\omega,\bar P)\rest 
[0,a_0))+\sum_{i<\omega} Th^m_{\bar k} ((\omega,\bar P)\rest [a_i,a_{i+1}))
=t+\sum_{i<\omega} s$.  

As $Th^n_r (\omega)$ is known, by \ref{2.8}, we can compute $Th^n_{\bar k}
(M,\bar P)$ from $s,t$. Now for any $t,s\in Th^n_{\bar k} (K^m_{\rm fin}),s\neq 
Th^n_{\bar k} (0,\bar P),\bar P\in \underline{P} (\varnothing)^m$, there is an 
$(\omega,\bar P)$ \st\ $Th^n_{\bar k} (\omega,\bar P)=t+\sum_{i<\omega} s$. 

As we know $Th^n_{\bar k} (K^m_{\rm fin})$ by \ref{3.3}, and can easily find whether 
$s\in Th^n_{\bar k} (K^m_{\rm fin})--Th^n_{\bar k} (\{0\})$, we finish.
\medskip

\begin{theorem}
\label{3.5} 
\begin{enumerate}
\item[(A)] From the monadic theory of $\lambda$ ($\lambda$ a cardinal) we can 
compute effectively the monadic theory of $K=\{\alpha:\alpha<\lambda^+\}$. 
\item[(B)] Moreover every monadic sentence which has model $\alpha<\lambda^+$, has 
a model $\beta<\lambda^\omega$.
\item[(C)] 
\begin{enumerate}
\item[(i)] For every $\alpha<\lambda^+$ there is a 
$\beta<\lambda^{\omega+1}+\lambda^\omega$ which has the same monadic theory
\item[(ii)] if $\mu\leqq\lambda$ and for every regular 
$\chi\leqq \lambda$ there is a $\chi'\leqq \mu$ \st\ 
$\chi,\chi'$ have the same monadic theory, then we can choose 
$\beta<\lambda^\omega \mu+\lambda^\omega$.\footnote{In fact, 
$\beta<M^{\omega+1}+M^\omega$.}
\item[(iii)] If we could always find $\chi<\mu$ then 
$\beta<\lambda^\omega \mu$, and if $\lambda=\omega, \beta
<\lambda^\omega+\lambda^\omega$.\footnote{In the first case 
$\beta<M$.}
\item[(iv)] Also, for every $\alpha<\lambda^+$, there are 
$n<\omega,\lambda_1,\ldots,\lambda_n\leqq \lambda$, \st\ 
the monadic theory of $\alpha$ is recursive in the monadic 
theories of $\lambda_1,\ldots, \lambda_n$, and 
$\lambda_i$ is a regular cardinal.
\end{enumerate}
\item[(D)] In general, the bounds in (B),(C) cannot be improved. 
\end{enumerate}
\end{theorem}

\par \noindent 
Remark:\ B\"uchi \cite{BS73} already proved (B),(C) for  % 2021-06-30 10:28 Bu73
$\lambda=\omega$ and (B) for $\lambda+\omega_1$.
\medskip

\par \noindent 
Proof:\ \begin{enumerate}
\item[(A)] Define $K_1=K_2=\{\alpha:\alpha\leqq\lambda\}$;
by \ref{3.2}(A)(i) and \ref{3.2}(B) % 2021-06-14 06:42 (A) 
    we can compute 
$Th^n_{\bar k} (K_i)$ for every $n,\bar k$ and 
$i=1,2$ (from the monadic theory of $\lambda$, of course).
Hence by \ref{2.5} we can compute 
$Th^n_{\bar k} (K')$ for every $n,\bar k$ where 
$K'=C\ell (K_1,K_2)$. Clearly every member of $K'$ 
is well-ordered and has cardinality $\leqq \lambda$. So up 
to isomorphism $K'\subseteqq K$. We should prove now only that 
equality holds. If not, let $\alpha$ by the first ordinal not 
in $K'$, and $\alpha<\lambda^+$. If $\alpha$ is a successor ordinal,
$\alpha -  % 2021-06-14 06:43 - 
    1\in K'; \text{ } 1,2\in K'$ % 2021-06-14 06:43 ' 
        hence $\alpha=(\alpha - 1)+1\in K'$, 
a contradiction. If $\alpha$ is a limit ordinal, its cofinality is 
$\leqq\lambda$. Let $\alpha=\sum_{i<i_0} \alpha_i,i_0\leqq 
\lambda,\alpha_i<\alpha$; then $i_0,\alpha_i\in K'$ so 
$\alpha\in K'$, a contradiction. 
\item[(B)] Let us first show that 
\begin{enumerate}
\item[(*)] For every $n,\bar k$ there is 
$q=q(n,\bar k)<\omega$ % 2021-06-14 06:44 )$ 
    \st\ 
if $\alpha,\beta<\lambda^+,\cf(\alpha)=\cf(\beta)$, and $\alpha,\beta$
are divisible by $\lambda^q$, then 
    $Th^n_{\bar{ k}} (\alpha)=Th^n_{\bar k}  % 2021-06-14 06:45 \bar{ }
(\beta)$.

For $n=0$ it is immediate, and we prove it for $n$. By the pigeon-hole
principle there are $1<\ell<p\leqq 2 |T(n,0,\bar k)|+1$ \st\ 
$Th^n_{\bar k} (\lambda^\ell)= Th^n_{\bar k} (\lambda^p)$. 
Clearly, 
$$
\lambda^{\ell+2}=\sum_{i<\lambda} (\lambda^{\ell+1}+\lambda^\ell).
$$
Hence
\[\begin{array}{ll}
Th^n_{\bar k} (\lambda^{\ell+2}) & = Th^n_{\bar k} 
[\sum_{i<\lambda} (\lambda^{\ell+1}+\lambda^\ell)]=\sum_{i<\lambda}
Th^n_{\bar k} (\lambda^{\ell+1}+\lambda^\ell) \\
\ & \quad \sum_{i<\lambda} [Th^n_{\bar k}(\lambda^{\ell+1})+
Th^n_{\bar k}(\lambda^\ell)]=\sum_{i<\lambda}[Th^n_{\bar k}(\lambda^p) \\
\ & \quad \sum_{i<\lambda} Th^n_{\bar k} (\lambda^\ell)=Th^n_{\bar k} 
(\sum_{i<\lambda} \lambda^\ell)=Th^n_{\bar k} (\lambda^{\ell+1}).
\end{array}\]

Hence we prove by induction on $m,\ell<m<\omega$ that 
$Th^n_{\bar k} (\lambda^m)=Th^n_{\bar k} (\lambda^{\ell+1})$; 
choose $q=q(n,\bar k)=\ell+1$. Let $\alpha,\beta<\lambda^+$ be divisible 
by $\lambda^q$ and have the same cofinality, and we shall prove $Th^n_{\bar k}
(\alpha)=Th^n_{\bar k} (\beta)$. Clearly it suffices to 
prove $Th^n_{\bar k} (\alpha)=Th^N_{\bar k} (\lambda^q \mu)$ where 
$\mu=\cf(\alpha)$. Let us prove it by induction on $\alpha$, and let 
$\alpha=\lambda^q \gamma$. If $\gamma=\gamma_1+1$, then for 
$\gamma_1=0$ i is trivial, and for $\gamma_1>0$
\[\begin{array}{ll}
Th^n_{\bar k} (\alpha) & = 
Th^n_{\bar k} (\lambda^q \gamma_1+\lambda^q)=
Th^n_k(\lambda^q \gamma_1)+Th^n_{\bar k}(\lambda^q) \\
\ & = Th^n_{\bar k} [\lambda^q \circ \cf (\lambda^q\gamma_1)]
+ Th^n_{\bar k} (\lambda^{q+2}) \\
\ & = Th^n_{\bar k} [\lambda^q \circ \cf (\lambda^q \gamma_1)+
\lambda^{q+2}]=Th^n_{\bar k} (\lambda^{q+2})=Th^n_{\bar k}
(\lambda^q\circ\lambda) \\
\ & = Th^n_{\bar k} [\lambda^q\circ \cf (\alpha)].
\end{array}\]
If $\gamma$ is a limit ordinal $\gamma=\sum_{i<\cf(\gamma)}
\gamma_i,\gamma_i<\gamma$ a successor,
\[\begin{array}{ll}
Th^n_{\bar k}(\alpha) & = 
Th^n_{\bar k}[\lambda^q(\sum_{i<\cf(\gamma)} \gamma_i)]=
Th^n_{\bar k}(\sum_{i<\cf(\gamma)} \lambda^q\gamma_i) \\
\ & = \sum_{i<\cf(\gamma)} Th^n_{\bar k}(\lambda^q \gamma_i) \\
\ & = \sum_{i<\cf(\gamma)} Th^n_{\bar k} 
[\lambda^q\circ \cf (\lambda^q\gamma_i)] \\
\ & \sum_{i<\cf(\gamma)} Th^n_k(\lambda^{q+1})=
\sum_{i<\cf(\gamma)} Th^n_{\bar k} (\lambda^q) \\
\ & = Th^n_{\bar k} [\lambda^q\circ \cf (\gamma)].
\end{array}\]
\end{enumerate}
So we have proved (*).
Let us prove (B). Let $\alpha<\lambda^+$ be a 
model of a sentence $\psi$. Choose by \ref{2.1}(A),(OR 3.2?) % 2021-06-14 06:47 xob
$n,  \bar k$
        % 2021-06-14 06:48 xob nosap- jo5el mahihu $ \beta $ 
        \st\ from $Th^n_{\bar k} (\beta)$ 
we know whether $\beta\models \psi$, and let 
$q=q (n,\bar k)$, and let $\alpha=\lambda^q \beta+
\gamma,\gamma<\lambda^q$. Then 
$$ 
Th^n_{\bar k}(\alpha)=Th^n_{\bar k}
[\lambda^q\circ\cf (\lambda^q \beta)+\gamma],
\mbox{ and } \lambda^q\circ \cf(\lambda^q\beta)+\gamma 
<\lambda^{q+2}.
$$
\item[(C)] Divide $\alpha$ by $\lambda^\omega$ so 
$\alpha=\lambda^\omega \alpha_1
+\alpha_2,\alpha_2<\lambda^\omega$. Let $\alpha'_1$ be 1
if $\alpha_1$ is a successor, and $\cf (\alpha_1)$ otherwise.  % 2021-06-14 06:29 cf ()
Then $\lambda^\omega \alpha_1,\lambda^\omega \alpha'_1$ 
are divisible by $\lambda^{q(n\bar k)}$ for every 
$n,\bar k$ and have equal cofinality. So by the proof of (B),
for every $n, \bar k, Th^n_{\bar k} (\lambda^\omega \alpha_1)=  % 2021-06-14 06:49 , 
Th^n_{\bar k} (\lambda^\omega \alpha'_1)$. Hence 
$\lambda^\omega \alpha_1+\alpha_2,\lambda^\omega \alpha'_1
+\alpha_2$ has the same monadic theory, and 
$\lambda^\omega \alpha'_1+\alpha_2<\lambda^\omega\lambda+
\lambda^\omega=\lambda^{\omega+1}+\lambda^\omega$. 
This proves (C)(i).

If $\chi'\leqq \mu$ has the same monadic theory as 
$\alpha'_1$ then $\lambda^\omega\alpha_1+\alpha_2,
\lambda^\omega \alpha'_1+\alpha_2$ and 
$\lambda^\omega \chi'+\alpha_2$ 
(which is $<\lambda^\omega \mu+\lambda^\omega$) 
have the same monadic theories. If 
$\chi'<\mu$ clearly $\lambda^\omega \chi'+
\alpha_2<\lambda^\omega \mu$.

If $\lambda=\omega$ then $\cf(\lambda) ^\omega \alpha_1)  % 2021-06-14 06:29 cf )
=\omega$ in any case, hence $\alpha=\omega^\omega \alpha_1
+\alpha_2$, and $\omega^\omega+\alpha_1<\omega^\omega+
\omega^\omega$ has the same monadic theory. Every 
$\alpha<\lambda^+$ we can uniquely represent as 
$$
\alpha=\lambda^\omega\alpha'+\lambda^n \alpha_n+
\ldots +\lambda^1 \alpha_1+\alpha_0; \alpha_i<\lambda.
$$
The monadic theory of $\alpha$ is recursive  in the monadic 
theories of $\lambda,\cf(\lambda)^\omega \alpha'),  % 2021-06-14 06:28 cf )
\alpha_n,\ldots,\alpha_0$. So we can prove inductively (C)(iv).
\item[(D)] Suppose $\lambda>\omega,\lambda$ is regular, and there is a 
sentence $\psi$ \st\ $\alpha\models \psi$ if $\alpha=\lambda$. 
Then there are sentences $\psi_n$ \st\ 
$\alpha\models \psi_n$ if and only if $\alpha=\lambda^n$, 
sentences $\varphi_n$ \st\ $\alpha\models \varphi_n$ if and only if 
$\alpha$ is divisible by $\lambda^n$, and sentence $\varphi$ \st\ 
$\alpha\models \varphi$ if $\cf( \alpha)  % 2021-06-14 06:28 cf () 
    =\lambda$. Then 
$\lambda^{\omega+1}$ is a model of $\{\varphi,\varphi_n:n<\omega\}$. 
If $\alpha$ is also a model of $\{\varphi,\varphi_n:n<\omega\}$ then 
$\lambda^n$ divides $\alpha$ for every $n$, hence $\lambda^\omega$ 
divides $\alpha$, so $\alpha=\lambda^\omega \beta$. If $\beta$ is a successor, 
$\cf (\alpha)=\omega$ but $\alpha\models \varphi$ so $\beta$ is a limit hence 
$\cf (\alpha)=\cf(\beta)$, so $\cf (\beta)=\lambda$, so   % 2021-06-14 06:28 cf()
$\beta \geqq \lambda$ hence $\alpha\geqq \lambda^\omega \circ \lambda
=\lambda^{\omega+1}$. Similarly $\lambda^{\omega+1}+\lambda^n$ is the smallest 
model of its monadic theory.
\end{enumerate}
\medskip

\begin{lemma}
\label{3.6} 
\begin{enumerate}
\item[(A)] In \ref{3.5}(A) it suffices to know the monadic theory of 
$\{\mu:\mu$ a regular cardinal $\leqq \lambda\}$. 
So if $\lambda$ is singular it suffices to know the monadic 
theory of $\{\alpha:\alpha<\lambda\}$. 
\item[(B)] For every sentence $\psi$,
\begin{enumerate}
\item[(1)] there is a sentence $\varphi$ 
(all in the monadic theory of order) \st\ $\alpha\models \varphi$ if and 
only if $\alpha$ is a limit and $\cf (\alpha)\models \psi$,
\item[(2)] there is a sentence characterizing the first 
ordinal which satisfies $\psi$ and 
\item[(3)] for every $n<\omega$ there is $\varphi_n$ \st\ 
$\alpha\models \varphi_n$ if and only if $\varphi$ is the $n^{\rm th}$
regular cardinal satisfying $\psi$. 
\end{enumerate}
\item[(C)] There are monadic sentences $\varphi_n$ \st\ $\alpha
\models \varphi_n$ if and only if $\alpha=\omega_n$. If $V=L$ 
there are monadic sentences  $\varphi^1_n$ \st\ $\alpha\models 
\varphi^1_n$ if and only if $\alpha$ is the $n^{\rm th}$ 
weakly compact cardinal. 
\end{enumerate}
\end{lemma}

\par \noindent 
Proof:\ 
\begin{enumerate}
\item[(A)] Immediate by \ref{3.5}(C)(iv). 
\item[(B)] \begin{enumerate}
\item Let $\varphi$ say that there is no last element, and for 
any unbounded $P$ there is an unbounded $Q\subseteqq P$ which satisfies 
$\psi$ (if $\cf ( \alpha)  % 2021-06-14 06:27 cf ()
    \models \neg \psi$ we can choose $Q$ as a set 
of order-type $\cf (\alpha)$;   % 2021-06-14 06:27 cf()
    so $\alpha\models \varphi$. 
If $\cf (\alpha)  % 2021-06-14 06:27 cf()
    \models \neg \psi$, let $P$ be a subset of $\alpha$ of 
order-type $\cf (\alpha)$; hence any unbounded $Q\subseteqq P$ has  % 2021-06-14 06:26 cf()
order-type $\cf (\alpha)$, so $\alpha\models \neg \varphi$).   % 2021-06-14 06:25 cf ()
\item Immediate. 
\item We use (1) and (2) to define $\varphi_n$ inductively. 
Let $\varphi_0$ say that $\alpha$ is the first ordinal whose cofinality 
satisfies $\psi$. Let $\varphi_{n+1}$ say that $\alpha$ is the first ordinal 
whose cofinality satisfies $\psi \wedge \neg \varphi_0 \wedge\ldots
\wedge \neg \varphi_n$. 
\end{enumerate}
\item[(C)] For $\varphi_n$ use (B)(3) for $\psi$ sating $\alpha$ is an 
infinite ordinal. For $\varphi^1_n$ use (B)(3) and Theorem \ref{0.1} 
(of Jensen). 
\end{enumerate}
\medskip

\section{The monadic theory of well-orderings}
If $a\in (M,\bar P)$ let 
$$
th(a,\bar P)=\{x\in X_i:a\in P_i\}\cup
\{x\notin X_i:a\notin P_i\}
$$
(so it is set of formulas). 

Let $D_\alpha$ denote the filter of (generated by) the closed 
unbounded subset of $\alpha,\cf ( \alpha)>\omega$.   % 2021-06-14 06:25 ()

\begin{lemma}
\label{4.1} 
If the cofinality of $\alpha>\omega$, then for every $\bar P\in 
\underline{P} (\alpha)^m$ there is a closed unbounded subset 
$J$ of $\alpha$ \st: for each $\beta<\alpha$, all the models 
$$
\{(\alpha,\bar P)\rest [\beta,\gamma):\gamma\in J,\cf(\gamma)
=\omega,\gamma>\beta\}
$$
have the same monadic theory. 
\end{lemma}  

\par \noindent 
Remark:\ 
B\"uchi \cite[6.1,p.110]{BS73} proved Lemma \ref{4.1} % 2021-06-30 10:29 Bu73
for $\alpha=\omega_1$, by a different method. 
\medskip

\par \noindent 
Proof: \ For every $n,\bar k$ there is, by \ref{1.1}, \ref{3.1}  % 2021-06-14 06:38 rewax
a homogeneous unbounded $I^n_{\bar k} \subseteqq \alpha$, by the coloring 
$f^n_{\bar k}$ of $(\alpha,\bar P)$, so there is $t^n_{\bar k}$ \st\ for 
every $\beta<\gamma\in I^n_{\bar k},\text{ } Th^n_{\bar k} ((\alpha,\bar P)\rest 
[\beta,\gamma))=t^n_{\bar k}$. Let $J^n_{\bar k}$ be the set of accumulation 
points of $I^n_{\bar k}$, and $J=\bigcap_{n,\bar k} J^n_{\bar k}$. 
Clearly $J$ is a closed and unbounded subset of $\alpha$. 

Let $\beta<\alpha$, and $\beta^n_{\bar k}$ be the first ordinal 
$>\beta$ in $I^n_{\bar k}$. Then for any $\gamma\in J,\gamma>\beta,\cf
(\gamma)=\omega$, and for every $n,\bar k$ we can find $\gamma_\ell\in 
I^n_{\bar k},\gamma_\ell<\gamma_{\ell+ 1 }, % 2021-06-14 06:40 J??},
    {\rm lim}_{\ell\rightarrow \omega}
\gamma_\ell= \gamma$ and $\gamma_0=\beta^n_{\bar k}$. Therefore 
\[\begin{array}{ll}
Th^n_{\bar k}((\alpha,\bar P)\rest [\beta,\gamma)) & = Th^n_{\bar k}
((\alpha,\bar P)\rest[\beta,\beta^n_{\bar k}))+\sum_{\ell<\omega}
Th^n_{\bar k}((\alpha,\bar P)\rest [\gamma_\ell,\gamma_{\ell+1})) \\
\ & = Th^n_{\bar k}((\alpha,\bar P)\rest [\beta,\beta^n_{\bar k}))
+\sum_{\ell<\omega} t^n_k.
\end{array}\]
So, $Th^n_{\bar k} ((\alpha,\bar P)\rest [\beta,\gamma)$ does not depend 
on the particular $\gamma$.
\medskip

\begin{definition}
\label{4.1a}
$A Th^n_{\bar k} (\beta,(\alpha,\bar P))$ for $\beta<\alpha,\alpha$ 
a limit ordinal of cofinality $>\omega$ is $Th^n_{\bar k}((\alpha,\bar P)
\rest [\beta,\gamma))$ for every $\gamma\in J,\gamma>\beta,\cf(\gamma)=\omega$;
where $J$ is from Lemma \ref{4.1}.
\end{definition}

\par \noindent 
Remark:\ As $D_\alpha$ is a filter, this definition does not depend on 
the choice of $J$.
\medskip

\begin{definition}
\label{4.2}
We define $W Th^n_{\bar k} (\alpha,\bar P)$:
\begin{enumerate}
\item if $\alpha$ is a successor or has cofinality $\omega$, it is 
$\varnothing$, 
\item otherwise we define it by induction on $n$:

for $n=0$: $W Th^n_{\bar k} (\alpha,\bar P)=
\{t:\{\beta<\alpha:th(\beta,\bar P)=t\}$ is a stationary 
subset of $\alpha\}$,
$$
\mbox{ for } n+1: \mbox{ let } W Th^{n+1}_{\bar k}
(\alpha,\bar P)=\{\langle S_1(\bar Q),
S_2(\bar Q)\rangle:\bar Q\in \underline{P} 
(\alpha)   % 2021-06-14 18:26 
    ^{{\bar k}(n+1)}\}
$$
where 
\[\begin{array}{ll}
S_1(\bar Q)=W Th^n_{\bar k} (\alpha,\bar P,\bar Q), \\
S_2 (\bar Q)=\{\langle t,s\rangle:\{\beta<\alpha:
W Th^n_{\bar k} ((\alpha,\bar P,\bar Q)\rest \beta)=
t,th(\beta,\bar P\conc \bar Q)=s\}
\end{array}\]
is a stationary subset of $\alpha\}$. 
\end{enumerate}
\end{definition}

\par \noindent 
Remark:\ Clearly, if we replace $(\alpha,\bar P)$ by a submodel 
whose universe is a closed  % 2021-06-14 18:30 class
    unbounded subset of $\alpha,
W Th^n_{\bar k} (\alpha,\bar P)$ will not change. 
Of course $W Th^n_k (M)$ is well defined for every 
well-ordered model.
\medskip

\begin{definition}
\label{4.3}
Let $\cf(\alpha)>\omega,M=  % 2021-06-14 18:25 )
    \alpha,\bar P)$ and we 
define the model $g^n_{\bar k} (M)=(\alpha,
g^n_{\bar k} (\bar P))$. 

Let 
$$
(g^n_{\bar k} (\bar P))_s=\{\beta<\alpha:s=
A Th^n_{\bar k} (\beta,M)\}
$$
and (when $m=\ell(\bar P)$)
$$
g^n_{\bar k} (\bar P)=\langle \ldots,
(g^n_k (\bar P))_s,\ldots
\rangle_{s\in T(n,m,\bar k)}.
$$
\end{definition}

\par \noindent 
Remark:\ \begin{enumerate}
\item In $g^n_{\bar k} (\bar P)$ we unjustly 
omit $\alpha$, but there will be no confusion.
\item Remember $T(n,m,\bar k)$ is the set of formally 
possible $Th^n_{\bar k} (M,\bar P),\ell(\bar P)=m$.
\end{enumerate}
\medskip

\begin{lemma}
\label{4.3a}
\begin{enumerate}
\item[(A)] $g^n_{\bar k} (\bar P)$ is a partition 
of $\alpha$. 
\item[(B)] $g^n_{\bar k} (\bar P\conc\bar Q)$ is a 
refinement of $g^n_{\bar k} (\bar P)$ and we can 
effectively correlate the parts. 
\item[(C)] $g^{n+1}_{\bar k} (\bar P)$ is a refinement 
of $g^n_{\bar k} (\bar P)$ and we can effectively correlate 
the parts. 
\item[(D)] The parallels of Lemma \ref{2.2a} 
for $Th,p Th$, hold for $W Th, p W Th$.
\end{enumerate}
\end{lemma}

\par \noindent 
Proof:\ Immediate. 
\medskip

\begin{theorem}
\label{4.3b} 
For every $n,m,\bar k$ we can effectively 
find $\bar r=\bar r_1 (n,m,\bar k)$ \st: 
If $\cf (\alpha^i)>\omega, M_i=(\alpha^i,
\bar P^i),\ell(\bar P^i)=m$ for 
$i=1,2$ and $A Th^n_{\bar k} (0,M_1)=
A Th^n_{\bar k} (0,M_2)$ and $W Th^n_{\bar r}
(g^n_{\bar k} (M_1  % 2021-06-14 18:33 i
    ))=W Th^n_{\bar r} 
(g^n_{\bar k} (M_ 2  % 2021-06-14 18:32  i % 2021-06-14 18:29 2
    ))$ then $Th^n_{\bar k} (M_1)
=Th^n_{\bar k} (M_2)$.
\end{theorem}

\par \noindent 
Proof:\ We prove by induction on $n$.   

For $n=0$, it is easy to check that 
$Th^n_{\bar k} (M_i)= A Th^n_{\bar k} 
(0,M_i)$ hence the theorem is trivial.

Suppose we have proved the theorem for $n$, and we shall
prove it for $n+1$. Suppose $\bar Q^1\in \underline{P}
(\alpha^1)^{\bar k(n+1)}$, and we shall find 
$\bar Q^2\in \underline{P} (\alpha^2)^{\bar k(n+1)}$ 
\st\ $Th^n_{\bar k}(\alpha^1,\bar P^1,\bar Q^1)=
Th^n_{\bar k} (\alpha^2,\bar P^2,\bar Q^2)$; 
be the symmetry in the hypothesis this is sufficient. 
Let $g^n_{\bar k} (\bar P^1\conc \bar Q^1)=\bar Q^{*1},
g^{n+1}_{\bar k} (\bar P^1)=\bar P^{*1},
g^{n+1}_{\bar k} (\bar P^2)=\bar P^{*2}$. Define  
$\bar r(n+1)=\ell(g^n_{\bar k} 
(\bar P^1\conc \bar Q^{-1}))=\ell(\bar Q^{*1})$ and 
$\bar r\rest (n+1)=r_1(n,m+\ell(\bar P^1),\bar k)$. 

By the assumptions and Definition \ref{4.2}, there is 
$\bar Q^{*2}\in \underline{P}
(\alpha^2)^{\bar k(n+1)}$ \st\ (for our 
$n,\bar r$ and $\alpha^2,\bar P^{*2}$; 
$\alpha^1,\bar P^{*1}$), $S_\ell(\bar Q^{*1})
=S_\ell(\bar Q^{*2})$ for $\ell=1,2$. 
(The notation is inaccurate, but should be clear.)
So, for $\ell=1$, we get $W Th^n_{\bar r}(\alpha^1,
\bar P^{*1},\bar Q^{*1})=W Th^n_{\bar r} (\alpha^2,
\bar P^{*2},\bar Q^{*2})$, and without loss of generality 
$0\in Q^{*1}_s\leftrightarrow 0\in Q^{*2}_s$. 
(From now on we can replace $\bar r$ by $\bar r\rest
(n+1)$.) So by Lemma \ref{4.2}, for 
$\ell=1,2,\bar Q^{*\ell}$ is a partition of 
$\alpha^\ell$ refining $\bar P^{*\ell}$, 
hence for every $\beta<\alpha^\ell$ there is a 
unique $s_\ell(\beta)$ \st\ 
$\beta\in Q^{*\ell}_{s_\ell(\beta)}$. 

Now, for $\ell=1,2$, choose a closed unbounded subset 
$J_\ell$ of $\alpha^\ell$ \st:
\begin{enumerate}
\item[(0)] every member of $J_\ell$ which is not an 
accumulation point of $J_\ell$, has cofinality $\omega$,
\item[(1)] for any $s$, if $Q^{*\ell}_s$ is not a 
stationary subset of $\alpha^\ell$ then $Q^*_s
\cap J_\ell=\varnothing$,
\item[(2)] if $\beta<\gamma<\alpha^\ell$; 
$\cf (\gamma)=\omega$ then 

$Th^{n+1}_{\bar k} ((\alpha^\ell,\bar P^\ell)\rest 
[\beta,\gamma))=A Th^{n+1}_{\bar k} (\beta,
(\alpha^\ell,\bar P^\ell))$ \quad 
(use Lemma \ref{4.1}),
\item[(3)] for every $\gamma\in J_\ell,
\cf(\gamma)=\omega$,
$$
Th^{n+1}_{\bar k} ((\alpha^\ell,\bar P^\ell)
\rest [0,\gamma)),A Th^{n+1}_{\bar k}
(0,(\alpha^\ell,\bar P^\ell)),
$$
\item[(4)] \quad \quad \quad if $Q^{*\ell}_s\cap 
J_\ell \neq \varnothing, \beta\in J_\ell$ \\
then there are $\gamma\in J_\ell,\gamma>\beta,s_\ell
(\gamma)=s$ \st\ $\{\xi\in J_\ell:\beta\leqq \xi\leqq \gamma\}$
is finite,
\item[(5)] for any $s,t$, if $\{\beta<\alpha^\ell:t=
W Th^n_{\bar r} ((\alpha^\ell,\bar Q^{*\ell})
\rest \beta),s=Th(\beta,\bar Q^{*\ell})\}$ is 
not a stationary subset of $J_\ell$, then it is disjoint 
to $J_\ell$.
\end{enumerate}
\medskip

\par \noindent 
Remark:\ Note that (5) just strengthens (1).

Now we define $\bar Q^2$ by parts. That is, for every 
$\beta<\gamma\in J_2\cup\{0\},\gamma$ is the 
successor of $\beta$ in $J_2$, we define $\bar Q^2
\rest [\beta,\gamma)$ \st
$$
s_2 (\beta)=Th^n_{\bar k}((\alpha^2,\bar P^2\conc
\bar Q^2)\rest [\beta,\gamma)).
$$
This is possible as by definition of $s_2(\beta),\beta\in 
Q^{*2}_{s_\ell(\beta)}$, hence
$$
s_i(\beta)\in A Th^{n+1}_{\bar k} (\beta,(\alpha^2,
\bar P^2)).
$$
We now prove 
\begin{enumerate}
\item[(*)] if $\beta<\gamma\in J_2\cup\{0\},
\cf(\gamma)=\omega$, then
$$
s_2(\beta)=Th^n_{\bar k}((\alpha^2,\bar P^2,
\bar Q^2)\rest [\beta,\gamma)).
$$
We prove it by induction on $\gamma$ for all $\beta$.
\begin{enumerate}
\item[(i)] By (0) the first $\gamma>\beta_1,\gamma\in 
J_2$ has cofinality $\omega$, and by the definition of 
$\bar Q^2 (*)$ is satisfied.
\item[(ii)] Let $\beta<\xi<\gamma,\xi\in J_2$, for no 
$\zeta\in J_2,\xi<\zeta<\gamma$, has cofinality $\omega$. 
Then by the induction hypothesis $Th^n_{\bar k}
((\alpha^2,\bar P_2,\bar Q^2)\rest [\beta,\xi))=s_2
(\beta)$ and
$$
Th^n_{\bar k}((\alpha^2,\bar P^2,\bar Q^2)\rest [\xi,\gamma))
=s(\xi).
$$
We should now show that $s_2(\beta)+s_2(\xi)=
s_2(\beta)$. So it suffice to find $\beta'<\xi'<
\gamma'\in J_1,s_1(\beta')=s_2(\beta),\cf ( \xi' ) % 2021-06-14 06:35 cf ()
=\omega=\cf (\gamma'),s_1(\xi')= s_2(\xi')$;   % 2021-06-14 06:35 cf ()
and by the definition of $\alpha^2,Q^{*1}_{s_2(\beta)}$ is a 
stationary subset of $\alpha^1$, hence for some $\beta'\in 
J_1,\beta'\in Q^{*1}_{s_2(\beta)}$ hence 
$s_2(\beta')=s_2(\beta)$. As $\xi\in J_2$,
$$
\{\zeta\in Q^{*2}_{s_2(\xi)}:W Th^n_{\bar r}
(\alpha^2,\bar P^{*2},\bar Q^{*2})=\varnothing\}
$$
is stationary, hence we can find 
$\xi'\in J_1,ch(\xi')=\omega,s_2(xi')=
s_2(\xi)$.
\item[(iii)] If $\gamma$ is an accumulation point 
of $J_2$ the proof is similar to that of (ii).
Choose $\xi_m,m<\omega,\beta<\xi_m<\xi_{m+1}<\gamma,
{\rm lim}_m \xi_m=\gamma,\cf(\xi_m)=\omega$, and 
$s_2(\xi_m)=s_2(\xi_{m+1})$ (use (4)). Then

\[\begin{array}{ll}
Th^n_{\bar k}(\alpha^2,\bar P^2,\bar Q^2)\rest [\beta,\gamma))
& = Th^n_{\bar k}((\alpha^2,\bar P^2,\bar Q^2)\rest 
[\beta,\xi^0)) \\
\ & +\sum_{m<\omega} Th^n_{\bar k}((\alpha^2,\bar P^2,
\bar Q^2)\rest [\xi_m,\xi_{m+1})) \\
\ & = s_2(\beta)+\sum_{m<\omega} s_2(\xi_0).
\end{array}\]
We should prove this sum is $s_2(\beta)$, and this is done 
as in (ii).
\item[(iv)] There are $\xi\in J,\beta<\xi<\gamma,
\gamma$ the successor of $\xi$ in $J_2$ and $\cf(\xi)>\omega$.
As before we can find $\beta'<\xi'<\gamma'\in J_1,s_1
(\beta')=s_2(\beta),W Th^n_{\bar r} ((\alpha^1,
\bar P^{1*})\rest \xi')=W Th^n_{\bar r} ((\alpha^2,
\bar P^{*2})\rest \xi),s_1(\xi')=s_2(\xi),\cf (\xi')  % 2021-06-14 06:34 cf ()
>\omega,\cf (\gamma')=\omega$. So clearly   % 2021-06-14 06:34 cf ()
$$
Th^n_{\bar k}((\alpha^2,\bar P^2,\bar Q^2)\rest 
[\xi,\gamma))=s_2(\xi)=s_2(\xi')=Th^n_{\bar k}
((\alpha^1,\bar P^2.\bar Q^1)\rest [\xi',\gamma')).
$$
Now also 
$$
Th^n_k((\alpha^2,\bar P^2,\bar Q^2)\rest [\beta,\xi))
=Th^n_{\bar k}((\alpha^1,\bar P^1,\bar Q^1)\rest
[\beta',\xi'))
$$
by the induction hypothesis on $n$ and on $\gamma$. 
\end{enumerate}

So we have proved (*) and $g^n_{\bar k}((\alpha^2,\bar P^2,\bar Q^2))=
(\alpha^2,\bar Q^{*2})$.
\end{enumerate}
Now by the induction hypothesis on $n$ it follows that 
$Th^n_{\bar k} (\alpha^1,\bar P^1,\bar Q^1)=
Th^n_{\bar k} (\alpha^2,\bar P^2,\bar Q^2)$.
\medskip

\begin{theorem}
\label{4.4}
If $\cf(\alpha)>\omega$, 
$$
t_1  = % 2021-06-14 18:34 -
    W Th^n_r(g^n_{\bar k} (\bar P)),t_2=A Th^n_{\bar k}
(0,  (\alpha,\bar P)) , \text{ }  % 2021-06-14 18:35 .
  \bar r=\bar r_1(n,\ell(\bar P),\bar k),
$$
then we can effectively compute $Th^n_{\bar k}(\alpha,\bar P)$ 
form $t_1,t_2$.
\end{theorem}

\par \noindent 
Proof: \ The proof is similar to that of \ref{4.3}.
\medskip

\begin{conclusion}
\label{4.5}
If $\lambda$ is a regular cardinal, and we know 
 % 2021-06-14 18:36 $a 
$A Th^n_{\bar k} (0,\lambda),  \text{ } W Th^n_{\bar{ r}}   % 2021-06-14 18:37 
    (\lambda) , \text{ } % 2021-06-14 18:38  % 2021-06-14 18:37 
(\bar r=r_1(n,0  % 2021-06-14 18:38 9
    ,\bar k))$, then we can compute 
$Th^n_{\bar k}(\lambda)$.
\end{conclusion}

\begin{lemma}
\label{4.6}
If $\lambda$ is a regular cardinal $>\omega,\bar r=r
(n,0,\bar k)$, then, letting 
$T_1=\{Th^n_{\bar r} (\mu):\omega<\mu<\lambda, 
\mu$ a regular cardinal$\}, T_2=\{Th^n_{\bar r}
(\alpha):\alpha<\lambda\}$, we can compute effectively 
$A Th^n_{\bar k} (0,\lambda)$ from $T_1$; and we can compute 
$T_1$ effectively from $T_2$. 
\end{lemma}

\par \noindent 
Proof:\ Let $T=\{t_1,\ldots,t_n\}$, and if $t_i=
Th^n_r(\mu)$ let $t'_i=Th^n_{\bar k} (\mu^q),
q=q(n,\bar k)$, % 2021-06-14 18:39  *
    (we can compute it effectively: see the 
proof of \ref{3.5}(B) for the definition of 
$q(n,\bar k))$ and let $t=t'_1+\ldots+t'_\ell$,
then 
$$
\sum_{m<\omega} t=t \omega = A Th^n_{\bar k} (0,\lambda).
\footnote{The second phrase is immediate by \ref{3.6}(3).}
$$
\medskip

\begin{conclusion}
\label{4.7}  % 2021-06-14 18:40 \label{ }
Let $\lambda$ be a regular cardinal. If the monadic theory of 
$\{\alpha:\alpha<\lambda\}$, and $\{W Th^n_{\bar k} 
(\lambda):n,\bar k\}$ are given then we can compute effectively 
the monadic theory of $\lambda$.
\end{conclusion}

\begin{lemma}
\label{4.8}
For a regular $\lambda,\{W Th^n(\lambda):n<\omega\}$ and the 
first-order theory of $M^\lambda=(\underline{P}
(\lambda) / D_\lambda,\cup,\cap,-,\varnothing,1,\ldots,
R^\lambda_t,\ldots)$ are recursive one in the other, where 
$R^\lambda_t (P,\bar Q)$ holds if and only if 

$\{\beta<\lambda:\beta\in P$, and for some 
$n,t=W Th^n((\lambda,\bar Q)\rest \beta)\}\neq \varnothing 
({\rm mod}\text{ } D_\lambda)$. % 2021-06-15 06:29 
\end{lemma}

\par \noindent 
Remark:\ Note that for every $t$ there is at most one possible $n$.
\medskip

\par \noindent 
Proof:\ Immediate, similar to the proof of Lemma \ref{2.1a}.
\medskip

\begin{conclusion}
\label{4.9}
If the monadic theory of $\{\alpha:\alpha<\lambda\}$ and the 
first-order theory of $M^\lambda$ are decidable, then so is 
the monadic theory of $\lambda$.

Using \ref{4.9} we can try to prove the decidability of the 
monadic theory of $\lambda$ by induction on $\lambda$.

For $\lambda=\omega$ we know it by \ref{3.4}.

For $\lambda=\omega_1$ the $R^{\omega_1}_t$'s are trivial, 
(because each $\beta<\omega_1$ is a successor or 
$\cf (\beta)=\aleph_0$, hence by Definition \ref{4.3}(1),   % 2021-06-14 06:33 cf ()
$R^{\aleph_1}_t (P,\bar Q)$ holds if and only if 
$t=\varnothing$).
So it suffices to prove the decidability of 
$(\underline{P}(\omega_1)/D_{\omega_1},\cap,\cup,-,
\varnothing,1)$. But by Ulam \cite{Ul30} this is an atomless 
Boolean algebra, so its theory is decidable. 
Hence we reprove the theorem of B\"uchi \cite{BS73}.  % 2021-06-30 10:29 {Bu73}.
\end{conclusion}

\begin{conclusion}
\label{4.10}
The monadic theory of $\omega_1$ is decidable. 

Now we can proceed to $\lambda=\omega_2$. 
Looking more closely at the proof for $\omega_1$, we see that 
$W Th^n_{\bar k} (\omega_1,\bar P)$ can be computed from the 
set of atoms in the Boolean algebra generated by the $P_i$
which are stationary subsets of $\omega_1$; and 
we can replace $\omega_1$ by any ordinal of cofinality 
$\omega_1$. So all the $R^{\omega_2}_t$ can be defined 
by the function $F/D_{\omega_2}$,
$$
F(I)=\{\alpha<\omega_2:\cf (\alpha)=\omega_1,\alpha \sminus    
I\cap \omega_2\notin D_\alpha\}.
$$
\end{conclusion}

\begin{conclusion}
\label{4.11} 
The first order theory of 
$$
M^{\omega_2}_1=(\underline{P} (\omega_2)/
D_{\omega_2},\cap,\cup,-,\varnothing,1,F/D_{\omega_2})
$$
is decidable if and only id the monadic theory of 
$\omega_2$ is decidable. 

Notice that $F(I\cup J)=F(I)\cup F(J)$, and that for 
$M^{\omega_2}_1$ to have a decidable theory, it suffices that it have 
elimination of quantifiers. For this it suffices 
\begin{enumerate}
\item[(*)] for any stationary $A\subseteqq \{\alpha<
\omega_2:\cf(\alpha)=\omega\}$ and $B,C$ \st\ 
$F(A)=B\cup C$ there are stationary $A',B',A=A'
\cup B', A'\cap B'=\varnothing, F(A') =A ({\rm mod} 
D_{\omega_2})$ and $F(B')=B({\rm mod} D_{\omega_2})$.
\end{enumerate}
\end{conclusion}

\begin{conjecture}  \label{4A}   % 2021-06-14 18:44 4(A) 
(*) is consistent with ZFC.
\end{conjecture}
  
\section{From orders to uniform orders}

An equivalence relation $E$ on an % 2021-06-14 18:44  
ordered set $N$ is 
{\em convex} if $xEy , \text{ } % 2021-06-14 18:45 x\ E y, 
x<z<y\in N$, implies  $ x Ey$,   % 2021-06-14 18:45 $x\ E z$, 
i.e., every equivalence class is convex. 
On $N/E=\{\alpha/E:a\in N\}$ a natural ordering is defined. 
If $J$ is a convex of a model $(M,\bar P)$ then 
$th(J,\bar P)$ is $\langle \ell,s_1,s_2\rangle$ \st\ if 
there is no last (first) element in $J,s_2=1  \text{ }
    (s_1=1)$, if 
$b$ is the last (first) element, $s_2=th(b,\bar P) \text{ }
(s_1=th
(b,\bar P))$ (for definition, see the beginning of Section 4) 
and $\ell={\rm min} (|J|,2)$.

\begin{definition}
\label{5.1}
\begin{enumerate}
\item $\kappa(M)$ is the first cardinal $\kappa$, \st\ 
neither $\kappa$ nor $\kappa^*$ is embeddable in $M$. 
\item $\kappa(K)$ is l.u.b. $\{\kappa(M):M\in K\}$. 
\end{enumerate}
\end{definition}

\begin{definition}
\label{5.2} 
We define for every $n,\bar k$, the class $U^n_{\bar k}$ 
and $U Th^n_{\bar k} ((M,\bar P))$ for $M\in U^n_{\bar k}$
\begin{enumerate}
\item $U^n_{\bar k}=\{(M,\bar P):M$ is dense order with no 
first nor last element and there are $t_0$ and a dense 
$I\subseteqq |M|$ \st\ for every $a<b\in I$:
$$
t_0=Th^n_{\bar k} ((M,\bar P)\rest (a,b)) \mbox{ and }
th(a,\bar P)=th(b,\bar P)\}.
$$
Now we define $U Th^n_{\bar k} (M, \bar P)$ be induction on $n$. 
        % 2021-06-15 06:30 .==> ,
\item $U Th^0_{\bar k} (M,\bar P)=Th^0_{\bar k} (M,\bar P)$.
\item $U Th^{n+1}_{\bar k} (M,\bar P)=\langle S_1,S_2,
{\rm com}\rangle$ where 
\begin{enumerate}
\item[(A)] $S_1=\{U Th^n_{\bar k} (M,\bar P,\bar Q):\bar Q\in 
        % 2021-06-15 06:31 ((  ))
\underline{P} (M)^{\bar k(n+1)}, (M,\bar P,\bar Q)\in U^n_{\bar k}\}$,
\item[(B)] Before we define $S_2$, we make some conventions:
\begin{enumerate}
\item[$(\alpha)$] $T_1(T_2)$   % 2021-06-15 06:31 [T_2]$ 
is the set of formally possible 
$th(J,\bar P^1),J\neq \varnothing$, and $\ell(\bar P^1)=
\ell(\bar P),(\ell(\bar P^1)=\ell(\bar P)+\bar k(n+1))$;
\item[$(\beta)$] $T_3=\{\langle\ell,s_1,t,s_2\rangle:
\langle\ell,s_1,s_2\rangle\in T_2,t\in T(n,\ell(\bar P)+
\bar k(n+1),\bar k)$ and $\ell=1$ if and only if $t$ is 
the ``theory'' of the empty model$\}$;
    \item[$(\gamma)$] If $\langle\ell,s_1,s_2\rangle\in 
T_1, \langle\ell',s'_1,t,s'_2\rangle \in T_3$ then 
$\langle\ell,s_1,s_2\rangle \leqq \langle\ell',s'_1,t, % 2021-06-15 06:33 ,
s'_2\rangle$ when: $\ell=\ell'$ and 
$s_1=1\Leftrightarrow s'_1=1, s_2=1\Leftrightarrow 
s'_2=1$ and $s_1\neq 1\rightarrow s_1\subseteqq s'_1,
s_2\neq 1\rightarrow s_2\subseteqq s'_2$;
    \item[$(\delta)$] At last let $\bar r=\bar r (n,\ell
(\bar P),\bar k)$ be from \ref{2.8}, 
$S_2=\{U Th^n_{\bar r} (M/E,\bar P^*,\bar Q^*):E$
 a non-trivial convex equivalence relation over 
$|M|,(M/E,\bar P^*,\bar Q^*)\in U  % 2021-06-15 06:34 Y
    ^n_{\bar r},\bar P^*=
\langle\ldots, P^*_t,\ldots\rangle_{t\in T_1}$, where 
$P^*_t=\{a/E:a\in |M|, th(a/E,\bar P)=t\}$ and 
$Q^*=\langle \ldots,Q^*_t,\ldots\rangle_{t\in T_3}$ 
is a partition of $|M| /E$ refining $\bar P^*$ and 
$\emptyset\neq Q^*_{t(1)} \subseteqq P^*_t$ implies 
$t(1)\leqq t\}$.
\end{enumerate}
\item[(C)] Com is + if $M$ is a complete order, and $\minus$
otherwise.
\end{enumerate}
\end{enumerate}
\end{definition}

\begin{lemma}
\label{5.1a} 
\begin{enumerate}
\item[(A)] From $Th^{n+2}_{\bar k} (M,\bar P)$ we can 
check whether $(M,\bar P)\in U^n_{\bar k}$ and compute 
$U Th^n_{\bar k} (M,\bar P)$.
\item[(B)] Also the parallel to \ref{2.2} holds.
\end{enumerate}
\end{lemma}

\begin{lemma}
\label{5.2a}
For every dense $N\in K, \|N\|>1,n,\bar k$, 
there is a convex submodel $M$ of $N$ which belongs to 
$U^n_{\bar k}, \|M\|>2$. 
\end{lemma}

\par \noindent 
Proof:\ By Theorem \ref{1.3}, and  \ref{2.3a} % 2021-06-15 06:43 jiniti xob \ref{4.1}(A).
\medskip

\begin{lemma}
\label{5.3}
Suppose $N$ is a dense order, $\kappa(N)\leqq \aleph_1$;
$I\subseteqq |N|$ is a dense subset, and for every 
$a<b\in I,t_0=Th^n_{\bar k} ((N,P)\rest 
[a,b))$. Then there is $t_1$ \st\ 
\begin{enumerate}
\item for every $a<b\in |N|,t_1=Th^n_{\bar k} 
((N,\bar P)\rest (a,b))$.
\item Moreover for every convex $J\subseteqq |N|$, 
with no first nor last element, 
$t_1=Th^n_{\bar k} ((N,\bar P)\rest J)$.
\end{enumerate}
\end{lemma}

\par \noindent 
Proof:\ Clearly it suffices to prove (2). 
Choose $a_0\in J\cap I$. Now define $a_n, % 2021-06-15 06:44 ,
    0<n<\omega$ 
\st\ $a_n\in J\cap I, a_n<a_{n+1}$ and 
$\{a_n:n<\omega\}$ is unbounded in $J$ 
(this is possible as $\kappa(N)\leqq \aleph_1$). 
Now define similarly, $a_n\in J\cap I, n$ 
a negative integer so that $a_{n-1}<a_n<a_0$ and 
$\{a_n:n$ is a negative integer$\}$ is unbounded 
from below in $J$. 

So, letting $Z$ be the integers, 
$$
Th^n_{\bar k} ((N,\bar P)\rest J)=
\sum_{n\in Z} Th^n_{\bar k} ((N,\bar P)\rest
[a_n,a_{n+1}))=\sum_{n\in Z} t_0 
\stackrel{def}{=} t_1.
$$
\medskip

\begin{theorem}
\label{5.4}
Let $M$ be an order, 
$\kappa(M)\leqq \aleph_1$. 
\begin{enumerate}
\item[(A)] Knowing $t$ and that $t=
U Th^n_{\bar k} (M,\bar P), (M,\bar P)\in 
U^n_{\bar k}$ we can effectively compute 
$F(t)=Th^n_{\bar k} (M,\bar P)$. 
\item[(B)] If $(M^i ,  % 2021-06-15 06:46 2--> i\
    bar P^i)\in U^n_{\bar k}$ 
for $i=1,2$, and $U Th^n_{\bar k} (M^1,\bar P^1)=
U Th^n_{\bar k} (M^2,\bar P^2)$ then 
$Th^n_{\bar k} (M^1,\bar P^1)=Th^n_{\bar k}
(M^2,\bar P^2)$. 
\end{enumerate}
\end{theorem}

\par \noindent 
Proof:\ Clearly (A) implies (B). 
So we prove (A) by induction on $n$.

For $n=0$ it is trivial.

Suppose we have proved the theorem for $n$,
and we shall prove it for $n+1$.

Let $U Th^{n+1}_{\bar k} (M,\bar P)=
\langle S_1,S_2,{\rm com}\rangle$. 
We should find 
$$
T=\{Th^n_{\bar k} (M,\bar P,\bar Q):
\bar Q\in \underline{P} (M)^{\bar k(n+1)}\}.
$$

If $t\in S_1$, then for some 
$\bar Q\in \underline{P} (M)^{\bar k(n+1)},
(M,\bar P,\bar Q)\in U^n_{\bar k}$ and 
$t=U Th^n_{\bar k} (M,\bar P,\bar Q)$, hence, 
by the induction hypothesis $F(t)=Th^n_{\bar k}
(M,\bar P,\bar Q)$, so $F(t)\in T$. We 
can conclude that $T'=\{F(t):T\in S_1\}\subseteqq 
T$. 

Now if $t^*\in S_2$, then there is a convex 
equivalence relation $E$ on $M$, \st\ 
$t^*=U Th^n_{\bar r} (M/E,\bar P^*,\bar Q^*)$ where the
conditions of $S_2$ are satisfied. If 
$Q^*_{\langle\ell,s_1,t,s_2\rangle} \neq \varnothing$, 
and $\ell>1$ implies $t\in T$ then we can define 
$\bar Q\in \underline{P} (M)$ \st\ for 
$a/E\in Q^*_{\langle\ell,s_1,t,s_2\rangle}$:
\begin{enumerate}
\item $U Th^n_{\bar k} ((M,\bar P,\bar Q)
\rest {\rm int} (a/E))=t$, 
\item $th(a/E,\bar Q)=\langle\ell,s_1,s_2\rangle$.
\end{enumerate}
\medskip

\par \noindent 
Remark:\ 
(1)\ can be done because by Lemma 
\ref{5.3}(2)  % 2021-06-15 07:40 \ref{4.3a}(2)??,
if ${\rm int} (a/E)\neq \varnothing$ then 
$$
Th^{n+1}_{\bar k} ((M,\bar P)\rest {\rm Int} 
(a/E))= Th^{n+1}_{\bar k} (M,\bar P)=T.
$$
Now clearly knowing $t^*$ we can compute 
$$
S(t^*)=\{t:Q^*_{\langle\ell,s_1,t,s_2\rangle}
\neq \varnothing,t \neq Th^n_{\bar k} (\varnothing), 
\mbox { for some } s_1,s_2\}
$$
where $\bar Q^*$ is an above. We can also compute 
$G(t)=Th^n_{\bar k} (M,\bar P,\bar Q)$. 
We know that $t\in S_2,S(t)\subseteqq T$, imply 
$G(t)\in T$. 

We know also that if

% 2021-06-15 07:38 \begin{enumerate} 
% 2021-06-15 07:38 \item[(i)] % 2021-06-15 07:37 $(i)
(i)    $t=Th^n_{\bar k}
((M,\bar P)\rest \{a\})$ for some $a\in M$,
and 

(ii) % 2021-06-15 07:38 \item[*ii)]   % 2021-06-15 07:37 (ii) 
    $t_1,t_2\in T$, 
% 2021-06-15 07:38 \end{enumerate} 
    
    then:
$\sum_{0\leqq n} (t_1+t)\in T$ and 
$\sum_{\stackrel{n<0}{n\in Z}} (t+t_2)\in 
T,t_1+t+t_2\in T$ and if com is $-,t_1+t_2\in 
T$ (where $Z$ is the set of integers) (we use the 
facts that $M$ is dense, $\kappa(M)\leqq \aleph_1$).

Now let $T^*$ be the minimal subset of 
$T(n,\ell(\bar P),\bar k)$ \st
\begin{enumerate}
\item[(a)] $T^*\supseteqq T'$,
\item[(b)] $t\in S_2,S(t)\subseteqq T^*$ imply
$G(t)\in T^*$,
\item[(c)] if $t_1,t_2\in T^*,t=Th^n_{\bar k}((M,P)\rest
\{a\})$ then $t_1+t+t_2\in T^*$; 
\item[(d)] if $t_2\in T^*,t_1=Th^n_{\bar k} 
((M<\bar P)\rest \{a\})$ for some $a\in M$ then 
$$
\sum_{0\leqq n<\omega} (t_2+t_1)\in T^*,
\sum_{\stackrel{n\leqq 0}{n\in Z}} (t_1+t_2)\in 
T^*;
$$
\item[(e)] if $t_1,t_2\in T_2$, com is $\minus$
then $t_1+t_2\in T^*$.
\end{enumerate}

It is easy to see that as $S_1,S_2$ are given and 
$T(n,\ell(\bar P),\bar k)$ is (hereditarily) finite and 
known, we can effectively compute $T^*$. So it suffices to 
prove that $T=T^*$ but as clearly $T^*\subseteqq T$ it 
suffices to prove:
$$
t\in T\Rightarrow t\in T^*.
$$

As $t\in T$, there is $\bar Q\in \underline{P} 
(M)^{\bar k(n+1)}$ \st\ $t=Th^n_{\bar k} (M,\bar P,\bar Q)$.
Define the equivalence relation $E$ on $M: a E b$ if and 
only if $a=b$ or, without loss of generality we assume that 
$a<b$, for every $a',b'\in M, a\leqq a'<b'\leqq b,Th^n_{\bar k}
(M,\bar P,\bar Q)\rest (a',b'))\in T^*$. It is easy to check that 
$E$ is a convex equivalence relation over $M$. Now we shall show that if 
$a\in M$, int$(A/E)\neq \varnothing$ then $Th^n_{\bar k}  % 2021-06-15 07:34 rewax
((M,\bar P,\bar Q)\rest {\rm int} (a/E))$ belongs to $T^*$. 
Choose $a_0\in a/E$, and then define $a_n,n\geqq 0$ \st\  % 2021-06-15 07:35 .==>/
$a_n<a_{n+1},\{a_n:0\leqq n<\omega\}$ is unbounded in 
${\rm int}(a/E)$. 
Without loss of generality $th(a_n,\bar P\conc\bar Q)=
s_0$ for every $n>0$. 
Hence 
\[\begin{array}{ll}
Th^n_{\bar k} ((M,\bar P\conc\bar Q)\rest 
\{x\in {\rm int} (a/E):a_0<x\}) \\
=\sum_{0\leqq n<\omega} [Th^n_{\bar k} 
((M,\bar P,\bar Q)\rest (a_n,a_{n+1}))+
Th^n_{\bar k} ((M,\bar P,\bar Q)\rest \{a_{n+1}\})].
\end{array}\]

By the definition of $E, Th^n_{\bar k} ((M,\bar P,\bar Q)\rest 
(a_n,a_{n+1}))\in T^*$, hence by (d),
$$
Th^n_{\bar k} ((M,\bar P,\bar Q)\rest 
\{x\in {\rm int} (a/E):a_0<x\})\in T^*.
$$
Similarly, 
$$
Th^n_{\bar k} ((M,\bar P,\bar Q)\rest \{x\in 
{\rm int} (a/E):x<a_0\})\in T^*.
$$
So by (c), 
$$
Th^n_{\bar k} (M,\bar P,\bar Q)\rest {\rm int} (a/E))
\in T^*.
$$

Similarly, by (c),(e) in 
$M/E$ there are no two successive elements, so 
$M/E$ is a dense order.

Define $\bar P^*=\langle\ldots,P_{\langle\ell,s_1,s_2\rangle},
\ldots\rangle, \bar Q^*=\langle\ldots, Q^*_{\langle\ell,s_1,t,s_2\rangle},
\ldots\rangle$ \st\ 
\begin{enumerate}
\item $a/E\in P_{\langle\ell,s_1,s_2\rangle}$ if and only if 
$th(a/E,\bar P)=\langle \ell,s_1,s_2\rangle$,
\item $a/E\in Q^*_{\langle\ell,s_1,t,s_2\rangle}$ if and 
only id $Th^n_{\bar k} ((M,\bar P,\bar Q)\rest {\rm int} 
(a/E))=t$; and $th(a/E,\bar P\conc \bar Q)=
\langle\ell,s_1,s_2\rangle$. 
\end{enumerate}

By Lemma \ref{5.2}, $(M/E,\bar P^*,\bar Q^*)$ either has 
only one element or it has an interval 
$(a/E,b/E)\neq \varnothing$ \st\ 
$(M/E,\bar P^*,\bar Q^*)\rest (a/E,b/E)\in U^n_{\bar r}$.

Now we prove $a E b$ and so show that this case does not occur 
and $E$ has one equivalence relation, hence $Th^n_{\bar k}
(M,\bar P,\bar Q)\in T^*$ and so we shall finish.

Let $a\leqq a'<b'\leqq b$, then let 
\[\begin{array}{ll}
J_2=\{c\in M:a'/E<c/E<b'/E\}, \\
J_1=\{c\in M:a'<c\in {\rm int}( a'/E) % 2021-06-15 06:51  (   ??
    \}, \\
J_3=\{c\in M:b'>c\in {\rm int} (b'/E)\}.
\end{array}\]

By (b), $Th^n_{\bar k} ((M,\bar P,\bar Q)\rest 
J_2)\in T^*$; by (d) 
$Th^n_{\bar k} ((M,\bar P,\bar Q)\rest J_i)\in T^*$ 
for $i=1,3$. Hence by (c) and (e) $Th^n_{\bar k}
((M,\bar P,\bar Q)\rest (a',b'))\in T^*$. So 
$a E b$, and we finish. 
\medskip
 
\begin{theorem}
\label{5.5}
\begin{enumerate}
\item[(A)] If $\kappa(K)\leqq \aleph_1$, 
and for every $M\in K$, there is 
$N\in K\cap U^{n+1}$ extending $M$, then 
from $U Th^{n+1}_{\bar k} (K)=
\{U Th^{n+1}_{\bar k} (M):M\in K\cap 
U^{n+1}_{\bar k}\}$, we can compute 
$Th^n_{\bar k}(K)$. Hence if 
$U Th^n (K)$ is recursive in $n$, 
then the monadic theory of $K$ is decidable. 
    \item[(B)] Suppose $\kappa(K)\leqq \aleph_1,K$ 
is closed under $M+N,\sum_{n<\omega} M,
\sum_{\stackrel{n\in Z}{n\leqq 0}} M_n,  
\sum_{i\in Q} M_i$ are convex submodels and 
division by convex equivalence relations. 
    Then from $U Th^n_{\bar r} (K) \text{ } % 2021-06-15 06:51 
    (\bar r=r
(n,0,\bar k))$ we can compute $Th^n_{\bar k}(K)$. 
Hence if $U Th^n (K)$ is recursive in $n$, then 
the monadic theory of $K$ is decidable. 
\end{enumerate}
\end{theorem}

\par \noindent 
Proof: \begin{enumerate}
\item[(A)] Immediate. 
\item[(B)] Essentially the same as the 
proof of 5.4.
\end{enumerate}
\medskip

\par \noindent
Remark:\ Of course there are other versions of (B), 
e.g., for a class of complete orders. 
\medskip

\section{Applications of Section 5 to dense orders}

\begin{definition}
\label{6.1}
$K_S$ is the class of orders $M$ \st\ no 
submodel of $M$ is isomorphic to $\omega_1$ or $\omega^*_1$    
or an uncountable subset of the reals\footnote{Those are the Specker 
orders; we get them from Aronszajn trees.}
\end{definition}

\begin{lemma}
\label{6.1a}
\begin{enumerate}
\item[(A)] $K_S$ satisfies the hypothesis of \ref{5.5}(B). 
Also no member of $K_S$ is complete, except the finite ones. 
\item[(B)] $K_S$ has uncountable members, but $M\in K_S$ 
implies $\|M\|\leqq \aleph_1$. 
\end{enumerate}
\end{lemma}

\par \noindent 
Proof:\ \begin{enumerate}
\item[(A)] Immediate. 
\item[(B)] The Specker orders. See e.g., 
\cite{Je71}\footnote{There is some overlapping between $S_1$ 
and $S_2$.} for existence. 
\end{enumerate}
\medskip

\begin{theorem}
\label{6.2} 
\begin{enumerate}
\item[(A)] The monadic theory of $K_S$ is decidable. 
\item[(B)] All dense order from $K_S$, with no first nor 
last element, have the same monadic theory.
\end{enumerate}
\end{theorem}

\par \noindent 
Proof:\ We shall show that for $(M,\bar P)\in U^0(K),
\bar P$ a partition, $p U Th^1 (M,\bar P)$ can be computed from 
$p U Th^0 (M,\bar P)$ (hence the former uniquely determine the 
latter). Then by the parallel to Lemma \ref{2.2a}, % 2021-06-15 07:42 , clause 
clause (B) follows immediately and 
(A) follows by\ref{5.5}(B). % 2021-06-15 07:43 (B)  % 2021-06-14 08:06 xob - lo5 barur \ref{5.6}.

So let $t=p U Th^0 (M,\bar P)$ be given; that is, we 
know that $\bar P$ is a partition of $M$ to dense or empty
subsets, $M\in U^0$, hence $M$ is dense with no first and 
no last element, $M\in K$, and we know $\{i:P_i\neq 
\varnothing\}$. So without loss of generality. 
$P_i\neq \varnothing$ for every $i$ and also 
$M\neq \varnothing,P_i$ is dense. Let $p Th^1
(M,\bar P)=\langle S_1,S_2,{\rm com}\rangle$, so we should 
compute com, $S_1,S_2$. 
\medskip

\par \noindent 
Part (1) com: \ As $M\in K$, and as clearly the rational order is 
embeddable in $M,M$ cannot be complete. 
\medskip

\par \noindent 
Part (2) $S_1$:\  It suffices to prove that any dense subset $P$ 
of $M$ can be split into two disjoint dense subsets of $M$.

So we shall prove more.
\begin{enumerate}
\item[(*)] If $M$ is a dense order, $I\subseteqq |M|$ is a dense subset, 
{\em then} we can partition $I$ to two dense subsets of $M$. That is, 
there are $J_1,J_2,I=J_1 \cup J_2, J_1\cap J_2=\varnothing$ and 
$J_1,J_2$ are dense subsets of $M$. 

We define a equivalence relation $E$ on $I: a E b$ if, 
$a=b$ or there are $a_0<a,b<b_0$ and $a_0<a'<b'<b_0$ implies 
$|\{c\in:a'<c<b'\}|=|\{c\in I:a<c<b\}|$ (and they are infinite by 
assumption). Now for every $E$-equivalence class $a/E$ with more 
than one element, let $\lambda=|\{a\in I:b'<a<c'\}|$ for every 
$b'\leqq c'\in a/E$.
\end{enumerate}
\medskip

\par \noindent
Case I:\ $|a/E|=\lambda>0$. 

Then let $\{\langle b_i,c_i\rangle:i<\lambda\}$ 
be an enumeration of all pairs $\langle b,c\rangle$ \st\ 
$b,c\in a/E,b<c$. Define by induction on $i<\lambda,
a^1_i,a^2_i\in a/E$. If we have defined them for 
$j<i$, choose
$$
a^1_i\in \{d\in I:b_i<d<c_i\} \sminus     \{a^2_j:j<i\},
$$
$$
a^2_i\in \{d\in I: b_i<d<c_i\} \sminus     \{a^1_j:j\leqq i\}.
$$
By cardinality considerations this is possible. Define 
$J_1 (a/E)=\{a^1_i:i<\lambda\}$. 
\medskip

\par \noindent 
Case II:\ $\lambda<|a/E|$.

Then clearly $|a/E|=\lambda^+$, and we can partition 
$a/E$ into $\lambda^+$ convex subsets 
$A_i,i<\lambda^+$, each of power $\lambda$. So 
on each we can define $J_1(A_1)$ \st\ 
$J_1 (A_i),A_i \sminus     J_1(A_i)$ are dense subsets of $A_i$. 
Let $J_1(a/E)=\bigcup_{i<\lambda^+} J_1 (A_i)$. 
\medskip

\par \noindent 
Case III:\ $\lambda=0$, so $|a/E|=1$. 

Let $J_1 (a/E)=\varnothing$. Let $J_1=
\bigcup_{a\in I} J_1 (a/E), J_2=I \sminus     J_1$. 

It is easy to check that $J_1,J_2$ are the desired subsets. 
\medskip

\par \noindent 
Part (3)\ $S_2$: By (2) it suffices to find to possible 
$U Th^0(M/E ,  % 2021-06-15 07:44 
    \bar P^*)$, where $\bar P^*=\langle \ldots, 
P^*_{\langle\ell,s_1, s_2  % 2021-06-15 07:45 , \rangle
\rangle},\ldots\rangle, 
P^*_{\langle  {\ell} % 2021-06-15 07:46 {\ell} 
    ,s_1,s_2\rangle}=
\{a/E  :  % 2021-06-15 07:46 ;==> :
    th(a/E,\bar P)=\langle \ell,s_1,s_2\rangle\}$, and 
$(M/E,\bar P^*)\in U^0 (K)$; so $W_E=
\{\langle \ell,s_1,s_2\rangle:P^*_{\langle\ell,s_1,s_2\rangle}
\neq \varnothing\}$ contain all relevant information. 
Clearly $W_E\neq \varnothing$ and $\langle\ell,s_1,s_2\rangle
\in W_E\Rightarrow \ell>0$ and we can also discard the case 
$\langle \ell,s_1,s_2\rangle\in W_E \Rightarrow \ell=1$. 
Also if $\langle \ell,s_1,s_2\rangle \in W_E$, then 
$\langle \ell,s_1,s_2\rangle$ is formally possible.

Suppose $W$ satisfies all those conditions, and 
we shall find a suitable $E$ \st\ $W_E=W$. Let 
$W=\{\langle \ell^i,s^i_1,s^i_2\rangle:i<q<\omega\}$. 
Choose a $J\subseteqq |M|$, countably dense in itself, unbounded in 
$M$ from above and from below, \st\ each $P_j\cap J$ is a dense
subset of $J$, and for no $a\in |M| \sminus     J$ is there a first (last)
element in $\{b\in J:b>a\}  \text{ }(\{b\in J; b<a\})$.  % 2021-06-15 07:47 rewax
$J$ defines $2^{\aleph_0}$ Dedekind cuts, but as $M\in K$, 
only $\leqq \aleph_0$ of them are realized. Let 
$\{a_n:n<\omega\}$ be a set of representatives from those cuts 
(that is, for every $a\in |M| \sminus     J$ there is $n<\omega$ \st\ 
$[a,a_n]$ or $[a_n,a]$ is disjoint to $J$). Let $J=\{b_n:n<\omega)$. 
Now we define by induction on $n$ a set $H_n$ of convex disjoint 
subsets of $M$, \st: 
\begin{enumerate}
\item[(a)] $H_n\subseteqq H_{n+1}; H_n$ is finite. 
\item[(b)] If $I_1\neq I_2\in H_n$ then $I_1<I_2$ or 
$I_2<I_1$ and between them there are infinitely many 
members of $J$. 
\item[(c)] If $I\in H_n,I$ has no last element, then for every 
$a\in |M| \sminus     J,a>I$, there is $b\in J,I<b<a$, and also 
$J\cap I$ is unbounded in $I$. 
\item[(d)] The same holds for the converse order. 
\item[(e)] If $I_1<I_2\in H_n, i<q$ then there are 
$I\in H_{n+1}, th(I,\bar P)=\langle\ell^i,s^i_1,s^i_2
\rangle$.\footnote{Also, $I_1<I<I_2$, and 
$I_0\in H_n$ implies $th(I_0,\bar P)\in W$.}
\item[(f)] $a_n,b_n\in \bigcup \{I:I\in H_n\}$. 
\item[(g)] If $I\in H_n$ has a first (last) 
element then this element belongs to $J$. It 
is not hard to define the $H_n$'s. Clearly 
$\bigcup_n \bigcup_{I\in H_n} I=
|M|$. So define $E$ as follows: 
$$
a E b \mbox{ if and only if } a=b \mbox{ or for some } 
n<\omega, I\in H_n,a,b\in I.
$$
It is not hard to check that $W_E=W$. So we finish the 
proof. 
\end{enumerate}
Along similar lines we can prove
\medskip
    
\begin{theorem}
\label{6.3}
Suppose $M$ is a dense order with no first nor last 
elements, $M$ is a submodel of the reals, and for every 
perfect set $P$ of reals, $P\cap |M|$ is countable, or 
even $<2^{\aleph_0}$. Then the monadic theory of $M$ is 
the monadic theory of rationals. 
\end{theorem}

\par \noindent 
Remark 1: 
We can integrate the results of \ref{6.2}, \ref{6.3}.
Always some $M$ satisfies the hypothesis of \ref{6.3}. 
If $2^{\aleph_0}>\aleph_1$, any dense $M\subseteqq R, |M|< 
2^{\aleph_0}$, and if $2^{\aleph_0}=\aleph_1$, the existence 
can be proved. 
\medskip

\par \noindent
Remark 2:\ In \ref{6.3} we can demand less of 
$|M|$: For all countable, disjoint and dense sets 
$Y_1,\ldots Y_n (n<\omega)$ there is a perfect set 
$P$ of reals \st\ $Y_i$ is dense in $P$ for 
$1\leqq i\leqq n$ and $P\cap |M|$ is $< 2^{\aleph_0}$ 
(see Section 7 for definition). 

The proof of \ref{5.2a}  % 2021-06-15 07:55 xob \ref{5.5}  % 2021-06-14 08:16 \ref{5.6} 
is easily applied to the 
monadic theory of the reals. (We should only 
notice that $R$ is complete.)
\medskip

\begin{conclusion}
\label{6.4} 
If we can compute the $U Th^n (R)$ for 
$n<\omega$ then the monadic theory of the real 
order is decidable. 
\end{conclusion}

\par \noindent 
Remark:\ Similar conclusions hold if we 
add to the monadic quantifier (or replace it by) 
$(\exists^{<\aleph_1} X)$ (i.e., there is a countable 
$X$). Notice that if $E$ is a convex equivalence relation 
over $R$, then $\{a/E:|a/E|>1\}$ is countable.

Grzegorczk \cite{Gr51} asked whether the lattice of subsets 
of reals with the closure operation has a decidable theory. 
One of the corollaries of Rabin \cite{Ra69} is that the theory 
of the reals with quantification over closed sets, and 
quantification over $F_\sigma$ sets is decidable. 

By our methods we can easily prove
\medskip

\begin{theorem}
\label{6.5}
The reals, with quantifications over countable sets, 
has a decidable theory. (We can replace ``$X$ countable''
by ``$|X|<2^{\aleph_0}$'' or ``$(\forall P)$ 
($P$ closed nowhere dense $\rightarrow |P\cap X|< 2^{\aleph_0}$)").
        % 2021-06-15 07:56 '')

As every closed set is a closure of a countable set, 
this proves again the result of Rabin \cite{Ra69} concerning  
        % 2021-06-30 10:42 Ra64 % 2021-06-14 18:01 
Grzegorczk's question. We can also prove by our method 
Rabin's stronger results, but with more technical difficulties. 
\end{theorem}

\section{Undecidability of the monadic theory of the real order}
Our main theorem here is 

\begin{theorem}
\label{7.1}
\begin{enumerate}
\item[(A)](CH) The monadic theory of the real order is undecidable.
\item[(B)](CH) The monadic theory of order is undecidable. 
\end{enumerate}
\end{theorem}

\begin{theorem}
\label{7.2}
(CH) The monadic theory of $K_n=\{(R,Q_1,\ldots,Q_n):Q_i\subseteqq \real \}$,
where the set quantifier ranges over countable sets, $1\leqq n$, is 
undecidable.
(We can even restrict ourselves to sets of rationals.)

Let $2^{\leqq\omega}$ be the set of sequences of ones and zeros of 
length $\leqq\omega$; let $\leqq$ be a partial ordering of 
$2^{\leqq \omega}$ meaning that it is an initial segment, 
$\prec$ the lexicographic order. 
\end{theorem}

\begin{theorem}
\label{7.3}
\begin{enumerate}
\item[(A)]CH The monadic theory of $(2^{\leqq \omega},\leqq,\prec)$ 
is undecidable.
\item[(B)](CH) The monadic theory of $K_n=\{(2^{\leqq\omega},\leqq,
\prec,Q_1,\ldots,Q_n):Q_i\subseteqq 2^{\leqq\omega}\}$, where the set
quantifier ranges over sets, $1\leqq n$, is undecidable.
(We can even restrict ourselves to subsets of $2^{<\omega}$).
\end{enumerate}

Instead of the continuum hypothesis, we can assume only:
\begin{enumerate}
\item[(*)] \quad ``The union of $<2^{\aleph_0}$ sets of the 
first category in not $R$''.

This is a consequence of Martin's axiom (see \cite{Je03}) hence weaker   % 2021-06-14 15:47 
than CH, but also its negation is consistent, (see Hechler \cite{He73}
and Mathias \cite{Mat74}  % 2021-06-30 09:40 
and Solovay \cite{So70}). Aside from countable 
sets, we can use only a set constructible from any well-ordering of the 
reals. Remember that by Rabin \cite{Ra69} quantification over closed and 
$F_\sigma$ sets gives us still a decidable theory.
\end{enumerate}
\end{theorem}

\par \noindent Conjecture 7A:\ The monadic theory of $(2^{\leqq\omega},
\leqq  , \prec)$, where the set quantifier ranges over Borel sets only,
is decidable. 

This should be connected to the conjecture on Borel determinacy 
(see Davis \cite{Da64}, Martin \cite{Mr} and Paris % 2021-06-30 10:35 Mar70
\cite{Pa72}).\footnote{Meanwhile Martin \cite{Mr75} proved the Borel determinacy.}  % 2021-06-14 17:23 
This conjecture implies
\medskip

\par \noindent Conjecture 7B:\ The monadic theory of the reals, 
where the set quantifier ranges over Borel sets, is decidable 
(by Rabin \cite{Ra69}).
\medskip

\par \noindent 
Conjecture 7C:\ We can prove \ref{7.1}-\ref{7.3} in ZFC.

Theorems \ref{7.1}(A),(B),\ref{7.3}(A) answer well known problems 
(see e.g., B\"uchi \cite[p.38, Problem 1,2a,2b,4a]{BS73}. Theorem  % 2021-06-30 10:35 Bu73
\ref{7.3}(B) answers a question of Rabin and the author.

Unless mentioned otherwise, we shall use CH or (*).
\medskip

\par \noindent
Notation:\ $\real$ denotes the reals. A {\em prefect} 
set is a closed, nowhere dense set of reals, with no 
isolated points and at least  % 2021-06-15 07:59 s 
two points (this is 
a somewhat  deviant definition). We use $P$ to denote 
prefect sets. Let $x$ be an inner point of $P$ if 
$x\in P$, and for every $\epsilon>0, (x \minus \epsilon,x)
\cap P\neq \varnothing, (x,x+\epsilon)\cap P\neq \varnothing$.
    Let $D  \subseteq  % 2021-06-15 07:59 S
    \subseteq R$ be dense in $P$ if for every inner point 
$x<y$ of $P$, there is an inner $z\in P\cap D,x<z<y$. 
Note that if $D$ is dense in $P,P$ is the closure of 
$P\cap D$. Real intervals will be denoted by 
$(a,b)$ where $a<b$, or by $I$; $(a,b)$ is an 
interval of $P$ if in addition $a,b$ are inner points if $P$.
\medskip

\begin{lemma}
\label{7.4}
Let $J$ be an index-set, the $D_i \text{ }(i\in J)$ countable  % 2021-06-15 08:01 
dense subsets of $R$, and $D=\bigcup_{i\in J} D_i$; and 
for every $P,|D\cap P|<2^{\aleph_0}$. Then there is 
$Q\subseteq \real\sminus     D, Q=Q\{D_i:i\in J\}$, \st\ 
\begin{enumerate}
\item[(A)] if $P\cap D\subseteqq D_i \text{ }(i\in J)$ and 
$D_i$ is dense in $P$ ($P$ is, of course, prefect) 
then $|P\cap Q|<2^{\aleph_0}$;
\item[(B)] if for no (interval) $I$ of $P$, 
and $i\in J,P\cap D \cap I\subseteqq D_i$ but 
$D$ is dense in $P$ then $P\cap Q \neq \varnothing$.
\end{enumerate}
\end{lemma}

\par \noindent 
Proof:\ Let $\{P_\alpha:0<\alpha<2^{\aleph_0}\}$ be any 
enumeration of the perfect sets. We define 
$x_\alpha,\alpha<2^{\aleph_0}$ by induction on $\alpha$. 

For $\alpha=0,x_\alpha\in R$ is arbitrary. 

For any $\alpha>0$, if $P_\alpha$ does not 
satisfy the assumptions of (B) then let 
$x_\alpha=x_0$ and if $P_ \alpha $  % 2021-06-15 08:02 _\alpha 
    satisfies the assumptions of (B) 
let $x_\alpha\in P_\alpha \sminus    \bigcup \{
P_\beta:\beta<\alpha, (\exists i\in J) 
(P_\beta\cap D\subseteqq D_i$ and $D$ is dense in 
$P_\beta)\}=D$.

This is possible because for any $\beta,i$, if 
$P_\beta\cap D \subseteqq D_i,D$ is dense in $P_\beta,
P_\beta\cap P_\alpha$ is a closed nowhere dense subset 
of $P_\alpha$. As otherwise for some interval $I$ of 
$P_\alpha,P_\beta\cap P_\alpha$ is dense in $P_\alpha$, so 
by the closeness  % 2021-06-15 08:03  ?? closedness
of $P_\beta\cap P_\alpha,P_\beta
\cap P_\alpha\cap I=P_\alpha\cap I$; therefore 
$$
D_i\supseteqq P_\beta \cap D \supseteqq P_\alpha
\cap I\cap D,
$$
a contradiction of the assumption on $P_\alpha$. 
So by (*) and the hypothesis $|P_\alpha\cap D|
<2^{\aleph_0}$ there exists such $x_\alpha$.

Now let $Q=\{x_\alpha:\alpha<2^{\aleph_0}\}$. 
If $P$ satisfies the assumption of (A), then 
$P\in \{P_\alpha:0<\alpha<2^{\aleph_0}\}$. 
Hence for some $\alpha,P=P_\alpha$, hence 
$P\cap D \subseteqq \{x_\beta:\beta<\alpha\}$, 
so $|P\cap D|<2^{\aleph_0}$. If $P=P_\alpha$ 
satisfies the assumption of (B) then $x_\alpha\in 
P_\alpha,x_\alpha\in Q$, hence $P_\alpha\cap Q
\neq \varnothing$. So we have proved the lemma. 
\medskip

\begin{lemma}
\label{7.5}
There is a dense $D\subseteqq R$ and 
$\{D_i:i\in J\},|J|=2^{\aleph_0}$ \st
\begin{enumerate}
\item $|D\cap P|<2^{\aleph_0}$ for every perfect $P$. 
\item The $D_i$ are pairwise disjoint.
\item $D_i\subseteqq  D,  D_i $ is dense.   % 2021-06-15 08:03 ,D
\end{enumerate}
\end{lemma}

\par \noindent 
Proof:\ Let $\{P_\alpha:\alpha<2^{\aleph_0}\}$ 
enumerate the perfect subsets of $R$, and let 
$\{I_n:n<\omega\}$ enumerate the rational intervals 
of $R$, and if $\alpha=\delta+n$ ($n<\omega,\delta$ 
a limit ordinal) choose $x_\alpha\in I_n\sminus    
\bigcup_{\beta<\alpha} P_\beta\sminus    \{x_\beta:\beta<\alpha\}$ 
and let $D=\{x_\beta:\beta<2^{\aleph_0}),D_\alpha=
\{x_{\omega \alpha+n}:n<\omega\}$.
\medskip

\par \noindent 
Notation: \ $J$ will be an index set; 
$[J]^n=\{U:U\subseteqq J, |U|=n\}$, and if 
$D_i$ is defined for $i\in J$, let 
$D_U=\bigcup_{i\in U} D_i$. Subsets of 
$[J]^n$, i.e., symmetric 
$n$-place relations over $J$, are denoted by $S$;
and if we know $\{D_i:i\in J\},Q_S$ will by 
$Q\{D_U:U\in S\cup [J]^{n-1}\}$ from \ref{7.4}. 
\medskip

\begin{definition}
\label{7.1a}
Let $\varphi_n(X,D,Q,I^*)$ be the monadic formula saying 
\begin{enumerate}
\item[(A)] $X$ is a dense set in $I^*$ and $X\subseteqq D$.
\item[(B)] For every interval $I\subseteqq I^*$, and sets 
$Y_i,\dots , i=1,n+1$, if $Y_i\cap I\subseteqq X$ and the $Y_i$ are 
pairwise disjoint and each $Y_i$ is dense in $I$ then there is a 
perfect set $P,P\cap Q=\varnothing$, and each 
$Y_i\cap I$ is dense in $P$.
\end{enumerate}
\end{definition}

\par \noindent 
Remark:\ We can represent the interval $I_0$ as a convex set.
\medskip

\begin{lemma}
\label{7.6}
Let $D,\{D_i:i\in J\}$ be as in \ref{7.5}, 
$I^*$ an interval, $S\subseteqq ]J]^n, Q_S=
Q\{D_U:U\in S\cup [J]^{n-1}\}$ as in \ref{7.4}.
Then for any set $X\subseteqq R,\real \models \varphi_n 
[X,D,Q_S,I^*]$ if and only if 
\begin{enumerate}
\item[(A)] $X$ is dense in $I^*,X\subseteqq D$,
\item[(B)] for any interval $I\subseteqq I^*$ 
there is a subinterval $I_1$ and $U\in S\cup
[J]^{n-1}$ \st\ $X\cap I_1\subseteq D_U$.
\end{enumerate}
\end{lemma}

\par \noindent 
Proof:\ 
(I) Suppose $\real \models \varphi_n [X,D,Q_S,I^*]$. 
Then by (A) from Definition \ref{7.1a}, $X$ is dense in   
    % 2021-06-15 08:06 7.1==>71a
$I^*,X\subseteqq D$ so (A) from here is satisfies. To prove 
(B)let $I\subseteqq I^*$ be an interval, and suppose that 
for no subinterval $I_1$ of $I$ and for no $U\in S \cup 
[J]^{n-1}$, does $X\cap I_1\subseteqq D_U$ hold, and we 
shall get a contradiction. Now we define by induction on 
$\ell,1\leqq \ell\leqq n+1$, distinct $i(\ell)\in J$ and 
intervals $I^\ell,0\leqq \ell\leqq n$ so that 
$I^0=I,I^{\ell+1}\subseteqq I^\ell$, and 
$X\cap D_{i(\ell)}\cap I^\ell$ is dense in $I^\ell$.

If we succeed, in Definition \ref{7.1a}(B), choose 
$I^{n+1}$ as $I$, and $X\cap D_{i(\ell)}\cap I^{n+1}$ 
as $Y^\ell$. So necessarily by $\varphi_n$'s definition there 
is a perfect $P$ \st\ $X\cap D_{i(\ell)} \cap I^{\ell+1}$ 
is dense in $P$ for $\ell=1,n+1$, and $P\cap Q_S=\varnothing$.
But this contradicts Lemma \ref{7.4}(B) by the definition of 
$Q_S$. So for some $\ell<n+1$ we cannot find appropriate 
$i(\ell+1),I^{\ell+1}$. So if we let $Y=(X-\bigcup_{k\leqq \ell}
D_{i(k)})\cap I^\ell$, for no $I^+\subseteqq I^\ell$ and no 
$i\in J$ is $Y\cap D_i\cap I^+$ dense; i.e., for every 
$i\in J,Y\cap D_i$ is nowhere dense.

If $\ell=n$, but $\{i(1),\ldots, i(n)\}\notin S$ let 
$D_{i(n)}\cap X \cap I^\ell=Y^1_n\cup Y^1_{n+1}$, where 
$Y^1_n, Y^1_{n+1}$ are dense subsets of $I^\ell$, and 
$Y^1_k=X\cap D_{i(k)} \cap I^\ell$, and we get a contradiction 
as before.

If $Y$ is not dense in $I^\ell$, it is disjoint to some 
$I^+\subseteqq I^\ell$, $X\cap I^+$, so $X\cap I^+\subseteqq 
\bigcup_{k<\ell} D_{i(k)}$.  % 2021-06-15 08:07 K)}$. 
    So $U=\{i(0),\ldots, i(\ell) \in 
S\cup [J]^{n-1},X\cap I^+\subseteqq D_U$, 
contradicting an assumption we made in the beginning of the proof.
Hence $Y$ is dense in $I^\ell$.

As $(\forall  % 2021-06-15 08:07 _
    i \in J) Y\cap D_i$ is nowhere dense also for 
every finite $U  ssu \cup
    \subseteqq J,Y\cap D_U$ is nowhere dense. 
So we can chose inductively distinct $i_m\in J$ and distinct 
$x_m\in Y\cap D_{i_m}$ \st\ $\{x_{(n+3)m+k}:m<\omega\}$ are 
dense subsets of $I^\ell$, for $0\leqq k<n+2$. If we let 
$Y^2_k=\{x_{(n+3)m+k}:m<\omega\}$ for $k\leqq n+1$, by 
Definition \ref{7.1a} there is a perfect $P$, \st\ 
$Y^2_k$ is dense in $P,P\cap Q=\varnothing$, and we get 
contradiction by \ref{7.4}(B) and the choice of the 
$x_m$'s.

As all the ways give a contradiction, we finish one implication. 

(II)\ Now we want to prove that $\real \models \varphi_n [X,D,Q,I^*]$ 
assuming the other side. 

Clearly $X\subseteqq D$, and $X$ is dense in $I^*$ 
(by condition (A) of Lemma \ref{7.6}). So 
condition (A) in Definition \ref{7.1a} holds. For condition (B)
of that definition let $I\subseteqq I^*$ be an interval, 
$Y_k\cap I\subseteqq X,Y_k$ dense in $I$ for $k=1, \dots ,n+1$ and  % 2021-06-15 08:08 
$k\neq \ell\Rightarrow Y_k\cap Y_\ell=\varnothing$. We should find 
a perfect $P$ \st\ $P\cap Y_k$ is dense in $P$ and 
$P\cap Q=\varnothing$. We can choose a $U\in S\cup [J]^{n-1}$ 
and $I_1\subseteqq I$ so that $X\cap I_1\subseteqq D_U$ 
(by the hypothesis). Choose a perfect $P$ \st\ each $Y_k$ is 
dense in $P$. As $D$ is as in % 2021-06-15 08:09 Definition 
    \ref{7.4}, either case   % 2021-06-15 08:09 wither
gives $|P\cap D|<2^{\aleph_0}$. 
\begin{enumerate}
\item[(*)] Now we can find perfect $P_\alpha(\alpha<2^{\aleph_0})$
\st\ each $Y_k \text{ }( 1 % 2021-06-15 08:10 a text
    \leqq k\leqq n+1)$ is dense in $P_\alpha$ and 
$\alpha\neq\beta$ implies $P_\alpha\cap P_\beta\subseteqq 
\bigcup^{n+1}_{k=1} Y_k$. 
\end{enumerate}
\medskip

\par \noindent 
Proof of (*):\   For $\eta$ a finite sequence of ones and zeros 
$X_\eta$ will be a set of closed-open intervals and singletons 
with endpoints in $\bigcup^{n+1}_{k=1} Y_k $, % 2021-06-15 08:11 1$, 
    which are pairwise 
disjoint. We define $X_\eta$ by induction on $\ell(\eta)$. Let 
$X_{\langle  % 2021-06-15 08:11 ??
    \rangle}=\{[a,b)\}$, where $a,b\in Y_1$, and if 
$X_\eta$ is defined, for each interval $[a,b)\in X_\eta$, choose a 
decreasing \sq\ $x^a_i \text{ } (i<\omega)$ whose limit is $a$, and  % 2021-06-15 08:12 
$x^a_0<b$ and $x^a_i\in Y_k$ if and only if $\ell(\eta)=k
\text{ }{\rm mod} \text{ }  n+1,1\leqq k\leqq n+1$. Let, for $m=0,1$: 
    % 2021-06-15 08:13 , % 2021-06-15 08:13 text X 2
$$
X_{\eta\conc\langle m\rangle}=\{(x^1_{i+1},x^a_i):
\mbox{ for some } b, [a,b)\in X_\eta \mbox{ and }   
i=m  \text{ } {\rm mod} \text{ } 2\}
$$
$$
\cup\{\{a\}: \mbox{ for some } b, [a,b)\in X_\eta, \mbox{ or } 
\{a\}\in X_\eta\}.
$$
For $\eta$ a \sq\ of ones and zeros of length $\omega, P_\eta=
\bigcap_{\ell<\omega} (\bigcup X_{\eta \rest  n}).$ % 2021-06-15 08:15 \rest 

Because $|P\cap D|<2^{\aleph_0}$ for some $\alpha,P_\alpha
\cap D \subseteqq \bigcup^{n+1}_{k=1} Y_k$; so by \ref{7.4} 
(and the choice of $Q$'s), $|P_\alpha\cap Q_S|<2^{\aleph_0}$. We
can find $P^\beta_\alpha (\beta<2^{\aleph_0})$ \st\ each $Y_k$ is 
dense in $P^\beta_\alpha$ and $\beta\neq \gamma\Rightarrow 
P^\beta_\alpha\cap P^\gamma_\alpha\subseteqq \bigcup^{n+1}_{k=1}
Y_k$. So for some $\beta,P^\beta_\alpha\cap Q \subseteqq 
\bigcup^{n+1}_{k=1} Y_k \subseteqq D$, but $Q\subseteqq \real
\sminus     D$ hence $P^\beta_\alpha \cap Q=\varnothing$, and we finish. 
\medskip

\begin{definition}
\label{7.2a} 
Let $\psi_n (X,D,Q,I^*)$ be the monadic formula saying 
\begin{enumerate}
\item[(A)] $\varphi_n (X,D,Q,I^*)$,
\item[(B)] for any interval $I_1\subseteqq I^*$, if $Y$ is 
disjoint to $X$ and dense in $I_1$ then $\neg \varphi_n 
(X\cup Y,D,Q,I_1)$. 
\end{enumerate}
\end{definition}

\begin{lemma}
\label{7.7}
Let $D,J,D_i,S,Q_S$ be as in \ref{7.6}. 
Then for any $X\subseteqq \real,  \real\models \psi_n
[X,D,Q_S,I^*]$ if and only if 
\begin{enumerate}
\item[(A)] $X$ is dense in $I^*,X\subseteqq D$, 
\item[(B)] for any interval $I\subseteqq I^*$ 
there is a subinterval $I_1$ and $U\in S\cup
\{V\in [J]^{n-1}: (\forall i\in J) (V\cup 
\{i\}\notin S)\}$ \st\ $X\cap I_1=D_U\cap I_1$. 
\end{enumerate}
\end{lemma}

\par \noindent 
Proof:\ \begin{enumerate}
\item[(I)] Suppose $\real \models \psi_n [X,D,Q_S,I^*]$, 
then clearly condition (A) holds. 
For condition (B) let $I\subseteqq I^*$ be an interval. 
By Definition \ref{7.2a}(A), $\real \models \varphi_n   % 2021-06-14 08:17 7,2a
[X,D,Q_S,I^*]$, 
    hence by Lemma \ref{7.6}% 2021-06-15 08:19 (1)
(B), $I$ 
has a subinterval $I_0$ \st\ $X\cap I_0\subseteqq D_U$ 
where $U\in S\cup [J]^{n-1}$. If $(D_U \sminus     X)\cap I_0$ 
is somewhere dense, let it be dense in $I_1\subseteq I_0$,
and let $Y=(D_U \sminus     X)\cap I_1$, which gives us a contradiction 
to Definition \ref{7.2a}% 2021-06-15 08:19 (1)
(B). If $U\in [J]^{n-1}$, and for 
some $i\in J,V=U\cup\{i\}\in S$, we can get a similar 
contradiction by $Y=(D_V \sminus     X)\cap I_0$ in the interval 
$I_0$ (as $D_i\subseteq D_V \sminus     X,Y$ is dense). We can 
conclude that: $U\in S$ or $U\in [J]^{n-1}$ and 
$U\cup\{i\}\notin S$ for every $i\in J$ and that 
$(D_U \sminus     X)\cap I_0$ is nowhere dense. 
Hence for some $I_1\subseteq I_0,(D_U \sminus     X)\cap 
I_1=\varnothing$ hence $X\cap I_1=D_U \cap I_1$. 
\item[(II)] Now suppose that conditions (A),(B) 
hold; by Lemma \ref{7.6} it is easy to see that 
$\real \models \psi_n [X,D,Q_S,I^*]$. 
\end{enumerate}
\medskip

\begin{definition}
\label{7.3a}
Let $\chi_1 (D,Q,I^*)$ be the monadic formula saying:
\begin{enumerate}
\item[(A)] $D$ is dense in $I^*,I^*$ an interval;
\item[(B)] if $I\subseteq I^*,X,Y$ are dense in $I$ 
and 
$$ 
\real \models \psi_1 [X,D,Q,I]\wedge \psi_1 [Y,D,Q,I]
$$ 
then for some $I_1\subseteq I$, 
$$
X\cap Y\cap I=\varnothing \mbox{ or } 
X\cap I_1=Y\cap I_1.
$$
\end{enumerate}
\end{definition}

\begin{lemma}
\label{7.8}
\begin{enumerate}
\item[(A)] If $D,\{D_i:i\in J\}$, are as in \ref{7.5} 
then for any interval $I^*,\real \models \chi_1
[D,Q_J,I^*]$. 
\item[(B)] If $\real \models \chi_1 [D,Q,I^*]$ then we can 
find $I\subseteqq I^*$, and $X_i,i<\alpha_0$ \st\ 
\begin{enumerate}
\item each $X_i$ is a dense subset of $I$ and 
$\real \models \psi_1 [X_i,D,Q,I]$,
\item if $I_0 \subseteqq I$, and $X\subseteqq I_0$ 
is dense in $I_0$ and $\real \models \psi_1 [X,D,Q,I_0]$ 
then there are $i<\alpha$ and $I_1\subseteqq I_0$ 
\st\ $X\cap I_1 = X_i\cap I_1$.
\end{enumerate}
\item[(C)] In (B), $|\alpha_0|$ is uniquely defined 
by $D,Q,I$. 
\end{enumerate}
\end{lemma}  

\par \noindent 
Proof:\ \begin{enumerate}
\item[(A)] By \ref{7.7} it is immediate.
\item[(B)] Let $\{X_i:i<\alpha\}$ be a maximal family 
satisfying (1) and (2) for $I=I^*$. If for some interval 
$I$ there are no subintervals $I^1$ and dense $X^*\subseteqq X 
\cap I^1$ \st\ $(\forall i<\alpha_0)$ ($X_i\cap X^*$ is 
nowhere dense)\footnote{and 
$\real \models \psi_1 [X^*,D,Q,I^1].$} we are finished. 
Otherwise we can choose inductively on $n$ 
intervals $I^n\subseteqq I^*$ disjoint to 
$\bigcup_{\ell<n} I^\ell$ and $X^*_n\subseteqq 
X\cap I^n$ \st\ $(\forall i<\alpha_0), X_i\cap 
X^*_n$ is nowhere 
dense\footnote{and $\real \models \psi_1 [X^*,D,Q,I^n].$}, 
and \st\ $\bigcup_{n<\omega} I^n$ is dense in $I$. 
Then we could have defined $X_{\alpha_0}=\bigcup_{n<\omega}
 D^*_n$, a contradiction. 
\item[(C)] Easy. 
\end{enumerate}
\medskip

\begin{definition}
\label{7.4a}
Let $\chi^n (Q_1,D,Q,I^*)$ be the monadic  formula 
saying 
\begin{enumerate}
\item[(A)] $D$ is dense in $I^*$, which is an interval. 
\item[(B)] Suppose $I_0\subseteqq I^*, X_\ell\subseteqq 
I_0 (\ell<n)$ and $\real \models \bigwedge_{\ell<n} \psi_1 
(X_\ell,D,Q,I_0)$.Then there is $I_1\subseteqq I_0$ 
\st\ for all $I_2\subseteqq I_1$
$$
\real \models \psi_n (\bigcup_{\ell<n} X_\ell,D,Q_1,I_1) 
\equiv \psi_n (\bigcup_{\ell<n} X_\ell,D,Q_1,I_2).
$$
\end{enumerate}
\end{definition}

\begin{lemma}
\label{7.9}
If $D, \{D_i:i\in J\}$ are as in Lemma \ref{7.5},
$S\subseteqq [J]^n$ then for any interval $I^*,R
\models \chi^n [Q_S,D,Q,I^*]$. 
\end{lemma}

\par \noindent
Proof:\ Immediate. 
\medskip

\begin{theorem}
\label{7.10}
The set $A_r$ is recursive in the monadic theory 
of order; where $A_r=\{\theta:\theta$ 
is a first order sentence which has an $\omega$-model i.e., 
a model $M$ \st\ $(|M|,R_1)$ is isomorphic to 
$(\omega,x+1=y)\}$. 
\end{theorem}

\begin{conclusion}
\label{7.11}
True first order arithmetic is recursive in the monadic 
theory of order.
\end{conclusion}

\par \noindent 
Proof:\ It suffices to define for every first order sentence $\theta$, 
a monadic sentence $G(\theta)$ so that $\real \models G(\theta)$ if and 
only if $\theta$ has an $\omega$-models. 

By using Skolem-functions and then encoding them by relations, we 
can define effectively the sentence $G_1(\theta)$ \st\ $\theta$ has an 
$\omega$-model if and only if $G_1 (\theta)$ has an $\omega$-model and 
$$ 
G_1(\theta)=(\forall x_1,\ldots,x_{n(0)}) (\exists 
x_{n(0)+1},\ldots,x_{n(1)}) (\bigvee_i \bigwedge_j \theta_{i j}),
$$
$\theta_{i j}$ is an atomic, or a negation of an atomic, formula;
only the relations $R_0,\ldots, R_{n(2)}$ appear in it; 
$R_0$ is the equality; and $R_i$ has $m(i)$-places. 

Define (where $X,Y,D,Q$ are variables ranging over sets, 
$I$ is a variable ranging over intervals and $x,y$ are 
individual variables): 
\begin{enumerate}
\item[(0)] $G_2(X_k=X_\ell)=(\forall I^1\subseteqq I^*) 
(\exists I^2 \subseteqq I^1) (X_k \cap I^2=X_\ell \cap I^2)$,
\item[(1)] $G_2 [G_\ell (X_{k(1)},\ldots, 
X_{k(m(\ell))})] = (\exists Y) (Y\subseteqq D \sminus     D^*
wedge \bigwedge^{m(\ell)}_{i=1} \psi_2 
(X_{k(i)} \cup Y,D,Q^\ell_i,I^*)$ (for $\ell<0$),
\item[(2)] $G_2 (\theta)=(\forall X_1,\ldots, 
X_{n(0)}) (\exists X_{n(0)+1},\ldots, X_{n(1)})$ \\
$(\forall I^0 \subseteqq I) (\exists I^*\subseteqq I^0) 
[\bigwedge^{n(0)}_{\ell=1} \psi_1 (X_\ell,D,Q^*,I^*) \bigwedge 
\bigwedge^{n(0)}_{\ell=1} X_\ell \subseteqq D^*$
$$
\rightarrow \bigwedge^{n(1)}_{\ell=n(0)+1} X_\ell \subseteqq 
D^* \cap \bigwedge^{n(1)}_{\ell=n(0)+1} \psi_1 (X_\ell,D,Q^*,I^*)\wedge
\bigwedge_i \bigvee_j G_2 (\theta_{i j})].
$$
\item[(4)] Let $\chi^*$ be the conjunction of the following formulas:
\begin{enumerate}
\item[$(\alpha)$] $D,D^*$ are dense in $I,D^*\subseteqq D$,
\item[$(\beta)$] $\chi_1 (D,Q^*,I)$,
\item[$(\gamma)$] $\chi^2(Q^i_\ell,D,Q^*,I)$.
Let us denote 
$$
\tilde{R}_1 (X,Y,Q^1_1,Q^2_1,I')=(X\subseteqq D^* \wedge
Y\subseteqq D^*\wedge X \cap Y=\varnothing \wedge
$$  
$$
\psi_1 (X,D,Q^*,I') \wedge \psi_1 (Y,D,Q^* \wedge,I')\wedge
(\exists Z) [Z\subseteqq D \sminus     D^* \wedge \psi_1 (Z,D,Q^*,I')\wedge
$$
$$
\psi_2 (X\cup Z,D,Q^1_1,I')\wedge \psi_2 (Y\cup Z,D,Q^2_1,I')]
$$
and 
\item[$(\delta)$] $\psi_1 (X_0,D,Q^*,I) \wedge X_0
\subseteqq D^* \wedge (\forall Y) [\psi_1 (Y,D,Q^*,I) \wedge
Y \subseteqq D^*\rightarrow (\exists Y_1) \tilde{R}_1 (Y,Y_1)] 
\wedge (\forall I' \subseteqq I) (\forall Y) \neg \tilde{R}_1 
(Y,X_0,Q^1_1,Q^2_1,I') \wedge$
$$
(\forall Y_1 Y_2 Y_3) (\forall I^0 \subseteqq I) [\tilde{R}_1 
(Y_1,Y_2,Q^1_1,Q^2_1,I^0) \wedge \tilde{R}_1 (Y_1,Y_3,Q^1_1,Q^2_1,I^0)
$$
$$
\rightarrow (\forall I^1 \subseteqq I^0) (\exists I^2 
\subseteqq I^1) Y_2 \cap I^2=Y_3\cap I^2].
$$
\item[$(\epsilon)$] The formula saying that if $(\delta)$ 
holds when we replace $Q^1_1,Q^2_1$ by $\tilde{Q}^1_1,
\tilde{Q}^2_1$ resp. then 
$$
(\forall X) (\forall Y) (\forall I' \subseteqq I) [\tilde{R}_1 
(X,Y,Q^1_1,I') \rightarrow \tilde{R}_1 (X,Y,\tilde{Q}^1_1,
\tilde{Q}^2_1,I')].
$$
\end{enumerate}
\item[(5)] $G(\theta)=(\exists Q^*,D,D^*,X_0,\ldots,
Q^i_\ell,\ldots) (\forall I) [\chi^*\wedge G_3 
(\theta)]$.

Now we should prove only that $\theta$ has an $\omega$-model 
if and only if $\real \models G(\theta)$.
\begin{enumerate} 
\item[(I)]  Suppose $ M $ is an $ \omega $-model 2021-06-14 
if and only if $ \real\vdash G( \theta )$.  % 2021-06-14 08:01 
Let $J=\omega+\omega,D_i (i<\omega+\omega)$ be countable,
pairwise disjoint, dense subsets of $R$. Choose symmetric 
and reflexive relations $S^i_\ell$ on $\omega+\omega$ so
that 
$$
M\models R_\ell (x_1,\ldots,x_{k(\ell)})\Leftrightarrow 
(\exists y\in \omega+\omega) \bigwedge^{k(\ell)}_{i=1} 
\langle y,x_i\rangle \in S^i_\ell \wedge y\notin \omega).
$$
To prove $\real \models G(\theta)$, let $D=\bigcup_{i<\omega+\omega}
D_i,D^*=\bigcup_{i<\omega} D,Q^i_\ell=Q_{(S^i_\ell)}, X_0=D_0$, 
and $Q^*=Q_{\omega+\omega}$. Let $I$ be any interval. It is not 
hard to check that under those assignments $\real \models x^* \wedge
G_3 (\theta)$.
\item[(II)] Now suppose $\real \models G(\theta)$. Let 
$Q^*,D,D^*,X_0,Q^i_\ell$ be such that $\real \models 
(\forall I) (\chi^* \wedge G_3 (\theta))$. BY (4) 
$(\beta)$, clearly $\real \models (\forall I) \chi_1 (D,Q^*,I)$. 
Hence by Lemma \ref{7.8}(B) there are $I$ and $D_i,i<\alpha$ 
satisfying (1),(2),(3) from \ref{7.8}(B). As 
$\real \models (\forall I) (\chi^* \cap G_3 (\theta))$, then 
in particular $\real \models \chi^*\wedge G_3 (\theta)$. 
By (4)$(\delta)$, $\real \models \psi_1 (X_0,D,Q^*,I)$, so 
we can choose $D_0=X_0$. (See the proof of \ref{7.8}.) By 
(4)$(\delta)$ we can also assume that $\real \models \tilde{R}_1
(D_n,D_{n+1})$ for $n<\omega$. By (4)$(\epsilon)$ necessarily 
$D_i\subseteqq D^* \Leftrightarrow i<\omega$.   

Let $\{\bar j_\ell:\ell<\omega\}$ enumerate all sequences 
$j=\langle j(1),\ldots,j(n(0))\rangle$ of natural numbers. 
As $\real \models G_3 (\theta)$ for every $\bar j_\ell$ we 
can choose $X_i=D_{j_\ell(i)}$, and so there is an assignment 
$X_i\rightarrow D^{\ell,i}$ for $n(0)<i\leqq n(1)$ showing 
that $\real \models G_3 (\theta)$. So we can define by induction on 
$n<\omega$ intervals $I_n$ so that: $I_{n+1}\subseteqq I_n,
I_0\subseteqq I$, and for every $n(0)<i\leqq n(1)$ for some 
$j_n(i)<\alpha_0,D^{\ell,i} \cap I_{n+1}=D_{j_n(i)} \cap 
I_{n+1}$.

Now we define a model $M:|M|=\omega$, and $M\models 
R_\ell [j(1),\ldots,j(m(\ell))]\Leftrightarrow$ for some 
$n,\real \models (\exists Y) [Y\subseteqq D \sminus     D^*\bigwedge 
\bigwedge^{m(\ell)}_{i=1} \psi_2 
(D_{j(i)} \cap Y,D,Q^i,I^n_n)]$.

It is easy to check that $\real \models \theta$. 
\end{enumerate}
\end{enumerate}
\medskip

\par \noindent 
Remark:\ By some elaboration, we can add to the definition of 
$A_r$ also the demand 
\[\mbox{``}R_2 \mbox{ is a well-founded two-place relation''}\]
(also for uncountable structures). Thus, e.g., there are 
sentences $\theta_n$, \st\ MA implies: 
$R\models \theta$ if and only if $2^{\aleph_0}=\aleph_n$. 
\medskip

\begin{theorem}
\label{7.11a}
The set of first-order sentences which has a model, is 
recursive in the monadic theory of 
$\{(\real,Q):Q\subseteqq \real\}$ where the set-variables range 
over subsets of the rationals. 
\end{theorem}

\par \noindent 
Remark:\ Notice that a quantification over $P$ \st\ 
$D$ is dense in $P$ can be interpreted by a quantification 
over $P\cap D$, as the property ``$x$ in the closure of 
$X$'' is first-order. Hence $\varphi_n ,  \psi_n$ are, in   % 2021-06-15 08:21 , 
our restricted monadic theory.

By \ref{7.10},\ref{7.11}, Theorems \ref{7.1},\ref{7.2} 
and \ref{7.3} are in fact immediate. Theorem \ref{7.1}(B) can
also be proved by the following observation of Litman \cite{Li72},
which is similar to \ref{3.6}(B)(1):
\medskip

\begin{lemma}
\label{7.12}
The monadic theory of the real order is recursive in the 
monadic theory of order. 
\end{lemma}

\par \noindent 
Proof:\ For every monadic sentence $\theta$ let 
$G(\theta)$ be the monadic sentence saying:

``If the set $X$ is completely ordered, is dense and has 
no first nor last elements then some $Y\subseteqq X$ 
has those properties and in addition $(Y,<)\models \theta$.''

As every complete dense order contains a subset isomorphic 
to $\real$, and any complete dense order $\subseteqq \real$ with no 
first nor last element is isomorphic to $\real$, clearly 
$\real\models G(\theta)$ if and only if $\theta$ is satisfied 
by all orders so our results is immediate.
\medskip

\par \noindent 
Conjecture 7D:\ The monadic theory of $\real$ and the (pure) 
second-order theory of $2^{\aleph_0}$ are recursive in each 
other.\footnote{Gurevich proved it when $V-L$.}
\medskip

\par \noindent 
Conjecture 7E:\ The monadic theory of $\{\real,Q)  % 2021-06-14 17:47 ??
    :Q\subseteqq \real\}$ 
with the set-quantifiers ranging over subsets of the rationals; 
and the (pure) second-order theory of $\aleph_0$ are 
recursive in each other. Gurevich notes that if $V=L$ the 
intersection of 7D,E holds.
\medskip

\par \noindent 
Conjecture 7F:\ The monadic theory  of order and the (pure) 
second-order theory, are recursive in each other.

In conjectures 7D,E,F use (*) or CH if necessary.
\medskip

\par \noindent 
Conjecture 7G:\ If $D_\ell$ is a dense subset of $\real$,
and for every $P,|P\cap D_\ell|<2^{\aleph_0},$ for 
$\ell=1,2$ then $(\real,D_1),(\real,D_2)$ have the same monadic 
theory.\footnote{Gurevich disproved it.} 

% ??
\bibliographystyle{amsalpha}
\bibliography{shlhetal}

\end{document}